\documentclass[draft]{amsart}
\newtheorem{theorem}{Theorem}[section]
\newtheorem{lemma}[theorem]{Lemma}
\newtheorem{proposition}[theorem]{Proposition}
\newtheorem{corollary}[theorem]{Corollary}

\newtheorem{remark}[theorem]{Remark}

\newcommand{\bC}{\mathbb{C}}

\newcommand{\bQ}{\mathbb{Q}}
\newcommand{\bR}{\mathbb{R}}
\newcommand{\bZ}{\mathbb{Z}}

\newcommand{\GL}{\mathrm{GL}}
\newcommand{\M}{\mathrm{M}}

\newcommand{\disp}{\displaystyle}

\def\disp{\displaystyle}
 \numberwithin{equation}{section}
\begin{document}

\baselineskip=17pt

\title[Arithmetic differential equations and $E$-functions]{Arithmetic differential equations and $E$-functions}
\author{Said Manjra}
\address{
   Department of Mathematics\\
   University of Ottawa\\
   585 King Edward\\
   Ottawa, Ontario K1N 6N5, Canada}
\email{manjra@math.net} 
\thanks{Work  supported by CRSNG}
\maketitle
\begin{abstract} Let $K$ be a number field. We give
an arithmetic characterization at infinity of the differential
operator of $K[x,d/dx]$ with minimal degree in $x$ annihilating a
given $E$-function. Such an operator is called an
\emph{$E$-operator}.
\end{abstract}
\section{ Introduction}
Let $ K$ be a number field and let $V_{0}$ be the set of all finite
places $v$ of $K$. For each $v\in V_{0}$ above a prime number
$p=p(v)$, we normalize the corresponding $v$-adic absolute value so
that $|p|_{v}=p^{-1}$ and we put $\pi_{v}=p^{-1/(p-1)}$. We denote
by $K_v$ the $v$-adic completion of $K$. We also fix an embedding
$K\hookrightarrow {\bC}$. For a real number $r>0$, and a
differential operator $\phi\in K[x,d/dx]$, we denote by
$R_v(\phi,r)$ the generic radius of convergence, bounded above by
$r$, of a basis of solutions of $\phi$ in a neighborhood of a
$v$-adic generic point of absolute value $r$. Recall that the
\emph{Fourier transform} $\mathcal F$ is the $K$-automorphism of
$K[x,d/dx]$ which satisfies ${\mathcal F}(x)=d/dx$ and ${\mathcal
F}(d/dx)=-x$. A power series $g=\sum_{n\ge 0}a_nx^n\in K[[x]]$
(resp. $F=\sum_{n\ge 0}a_nx^n/n!\;\in K[[x]]$) is said to be a
\emph{$G$-function} (resp. an \emph{$E$-function}) if there exists a
positive constant $C$ such that for any index $n$, the coefficient
$a_n$ and its conjugates over ${\bQ}$ do not exceed $C^n$ in
absolute value, and if there exists a common denominator $d_m\ge 1$
for $a_0,\ldots,a_n$ which does not exceed $C^n$. Shudnovsky proved
in \cite{CC} that the minimal differential operator of $K[x,d/dx]$
annihilating a given $G$-function satisfies the \emph{Galochkin
condition} and such an operator is called later a
\emph{$G$-operator} \cite[IV]{A1}. E. Bombieri proved in 1982 that
the differential operator which have $G$-function solutions near
every regular singularity satisfies the condition $\prod_{v\in
V_0}R_v(\phi,1)\ne 0$ (called \emph{Bombieri's condition})
\cite[10]{Bo}. The equivalence between the condition of Galochkin
and that of Bombieri was established in 1989 by Y. Andr\'e
\cite[IV]{A1}. In 2000, the latter showed that the differential
operator of $K[x,d/dx]$ with minimal degree in $x$ annihilating an
$E$-function is the Fourier transform of a certain $G$-operator and
he called such operators \emph{$E$-operators} \cite[4]{A2}.
Recently, in joint work with Remmal \cite{MR}, we gave a local
$p$-adic characterization of the $E$-operators in the neighborhood
of $0$, which is a regular singularity. This result is given
 in term of the generic radius of convergence and provides
  an answer to a conjecture of Y. Andr\'e  \cite[4.7]{A2}.
In the present paper, we propose a local arithmetic characterization
of  the $E$-operators at infinity (Theorem 3.1), which is in general
an irregular singularity of such an operator. This result is the
analogue of the \emph{Local Bombieri} property for the $G$-operators
\cite[6]{CD}. In the proof of this result, we cannot avoid the case
of negative exponents as in \cite[6]{MR}. This requires the standard
Laplace transform instead of the formal one used in \cite[5]{MR}.

The importance of $E$-operators comes from the fact that if $y(x)$
is  an  \emph{arithmetic Gevrey series  of nonzero order} $s$ and is
a solution of  a linear differential equation with coefficients in
$K(x)$, then $y(x^{-s})$ is a solution  of  an $E$-operator  (cf.
\cite[6]{A2}).

This article is organized as follows:\\
In the seconde section, we start by giving some preliminaries which
will be needed later. In section 3, we state our main theorem
(Theorem 3.1), we give some key lemmas and we prove that the
conditions of Theorem 3.1 are necessary. The section 4  is devoted
to the Laplace transform $\mathcal{L}$;  in paragraph 4.1, we
summarize  main  formal properties of $\mathcal{L}$. In paragraph
4.2, we give some  arithmetic properties of $\mathcal{L}$. For a
given differential operator $\psi \in K[x,d/dx]$, we see in section
5, how we can determine the nature of solutions of  $d/dx.\psi$  at
$0$ from those of $\psi^*$ at the same point. Using the results of
sections 3, 4 and 5, we prove, in section 6, that the conditions of
Theorem 3.1 are sufficient.

\section{ Notations and preliminaries }
\subsection{2.1 Differential modules} \ \\
Let $\mathcal K$ be a commutative field equipped with a derivation
 $\partial$, let $K$ be  the constant field of $\partial$ in $\mathcal K$
and let $\mu$ be a positive integer. A differential $\mathcal
K$-module
 $\mathcal M$  is a free module of  rank $\mu$ over $\mathcal K$ equipped
with an $K$-endomorphism $\nabla$ of $\mathcal M$ which satisfies
 the condition $\nabla(am)=a\nabla(m)+\partial(a)m$ for any $m\in\mathcal M$ and
$a\in\mathcal
   K$. To each basis
$\{e_{i}\}$ of $\mathcal
   M$ over $\mathcal K$
  corresponds a matrix $G=(G_{ij})\in {\M}_\mu({\mathcal K})$
  satisfying
  $$\nabla(e_{i})=\sum_{j=1}^{\mu}G_{ij}e_{j}$$ called the
 matrix of $\partial$ with respect to the basis $\{e_{i}\}$
  (or simply the associated matrix of $\mathcal M$) and a
  differential system $\partial X=GX$ where $X$ denotes a column
  vector $\mu\times 1$  or $\mu\times \mu$ matrix.
 A change of  bases in $\mathcal M$ results in the existence of a matrix $Y\in
{\GL}_\mu({\mathcal
 K})$ such that  $Y[G]:=YGY^{-1}+\partial(Y)Y^{-1}$  is the associated matrix of
$\partial$
 in the new basis.
If $\phi= \disp\sum_{i=0}^{\mu} a_{i}\partial^{i}\in {\mathcal
K}[\partial]$ is a differential operator such that $ a_{\mu}\ne 0$,
one can associates  to it the differential ${\mathcal K}$-module
${\mathcal M}_{\phi}={\mathcal K}[\partial]/{\mathcal
K}[\partial]\phi$
 of rank $\mu$ which corresponds to a system $$
\partial X=A_{\phi}X\;\;\mbox{where}\;\; A_{\phi}:=
\left (
\begin{array}{c}
0\\0\\\vdots\\0\\-\frac{a_0}{a_{\mu}}
\end{array}
\begin{array}{c}
1\\0\\\vdots\\0\\\frac{a_1}{a_{\mu}}
\end{array}
\begin{array}{c}
0\\1\\\vdots\\0\\-\frac{a_2}{a_{\mu}}
\end{array}
\begin{array}{c}
\ldots\\\ldots\\\ddots\\ \\\ldots
\end{array}
\begin{array}{c}
0\\0\\\vdots\\1\\-\frac{a_{\mu-1}}{a_{\mu}}
\end{array}
\right)
$$
 is called the companion matrix of $\phi$. One associates to
$\phi$ the adjoint operator $\phi^{*}=
\disp\sum_{i=0}^{\mu}(-\partial)^{i}a_{i}$. One verifies that
$-^{T}A_{\phi}$ is associated to  ${\mathcal M}_{\phi^{*}}={\mathcal
K}[\partial]/{\mathcal K}[\partial]\phi^{*}$. More generally, $G$ is
associated to ${\mathcal M}_{\phi}={\mathcal K}[\partial]/{\mathcal
K}[\partial]\phi$ if and only if $-^{T}G$ is associated to
${\mathcal M}_{\phi^{*}}$. This comes from the fact that for any
$Y\in {\GL}_{\mu}({\mathcal K})$, one has
\begin{equation}
\label{S}
\begin{aligned}
-^{T}(Y[A_{\phi}]) &=
^{T}Y^{-1}(-^{T}A_{\phi})^{T}Y-^{T}Y^{-1}\partial(^{T}Y)\\
&=^{T}Y^{-1}(-^{T}A_{\phi})^{T}Y+\partial(^{T}Y^{-1})^{T}Y=(^{T}Y^{-1})[-^{T}A_{\phi}].
\end{aligned}
\end{equation}
\subsection{2.2 The Newton-Ramis polygon}\ \\ Let $\phi=\disp\sum_{i=0}^{\mu}
a_{i}(x)\Big(\disp\frac{d}{dx}\Big)^{i}= \disp\sum_{i=0}^{\mu}
\disp\sum_{j=0}^{\nu}a_{i,j}x^{j}\Big(\disp\frac{d}{dx}\Big)^{i}\in
K[x,\disp\frac{d}{dx}]$ be a differential  operator of rank $\mu$.
The  Newton polygon in the sense of Ramis of $\phi$, which we shall
denote by $NR(\phi)$, is the convex hull, in the plane $uv$, of the
horizontal half-lines $\{u\leq i,\;v=j-i\;|\;a_{i,j}\ne
 0\}$ (cf. \cite{Ra}).

With this definition, it is easy to check that
$NR(\overline{\phi})=NR(\phi)$ (where $\overline{\phi}$ denotes the
operator obtained from $\phi$ by the change of variable $x\to -x$).
Also, $NR(\phi)$ has non vertical side if and only if $a_{\mu}$ is a
monomial, in which case $\phi$ has non nonzero finite singularity.

The part of $NR(\phi)$ located in the half-plane $v\le ord_x(a_\mu)$
corresponds to the classic Newton polygon $N(\phi)$ of $\phi$. As
for the part of $NR(\phi)$ located in the half-plane $v\ge
deg(a_\mu)$, it corresponds, by translation, to the transform, by
the symmetry $(u,v)\longrightarrow (u,-v)$, of the Newton polygon
$N(\phi_\infty)$ of  the operator $\phi_\infty$ obtained from $\phi$
by the change of variable $x\to 1/x$ (cf. \cite[V]{Ma}). The slopes
of $N(\phi_\infty)$ are called the slopes of $\phi$ at infinity.

This implies that the non-vertical slopes of $NR(\phi)$ depend only
of $\mathcal{M}_\phi$, since $N(\phi)$ and $N(\phi_\infty)$ depend
only of $\mathcal{M}_\phi$ (cf. \cite[3.3.3]{VS}).

The polygon of $\mathcal{F}(\phi)$ may be obtained from $NR(\phi)$
by applying to it the transformation $(u,v)\to(u+v,-v)$ (cf.
\cite[V]{Ma}). This implies, in particular, that: \emph{$NR(\phi)$
has non nonzero finite slopes if and only if all slopes of
$NR({\mathcal F}(\phi))$ lie in $\{0,-1\}$.}

If  $A_1,\ldots,A_\ell$ are square matrices, we denote by
$\oplus_{1\le i\le \ell}A_{i}$ the block diagonal matrix
$$\oplus_{1\le i\le \ell}A_{i}= \left (
\begin{array}{c}
A_1\\ \ \\ \ \\ \
\end{array}
\begin{array}{c}
\ \\A_2 \\ \ \\ \
\end{array}
\begin{array}{c}
\ \\ \ \\ \ddots\\ \
\end{array}
\begin{array}{c}
\ \\ \ \\ \ \\ A_\ell
\end{array}
\right)
$$
with blocks $A_1,\ldots,A_\ell$ on the diagonal.
\subsection{2.3 Radius of convergence in neighborhood of
singularities}\ \\
Consider  the differential  field ${\mathcal K}=K(x)$ equipped with
the  derivation $\partial=d/dx$. Let $\phi$ be a differential
operator of rank $\mu$ such that the  slopes of $N(\phi)$ lie in
$\{0,1\}$ and let $G\in {\M}_\mu(K(x))$ be an associated matrix of
$\mathcal{M}_\phi$. The \emph{Turrittin-Levelt} decomposition states
that there exist  a finite extension $K'$ of $K$, a matrix
$Y_{0}(x)\in {\GL}_{\mu}(K'((x)))$, called a reduction matrix of
 $G$ (or simply of $\phi$ if $G=A_{\phi}$) at $0$, an upper triangular
matrix $C_{0}\in
  {\M}_{\mu}(K')$ and a diagonal  matrix $\Delta_{0}\in
  {\M}_{\mu}(K')$
commuting with $C_{0}$ such that
$Y_{0}(x)[G(x)]=\Delta_0/x^{2}+C_{0}/x$ \cite[3]{Le}. By base
change, we may assume that $C_0$ is in Jordan form.

One observes that the matrix
$Y_{0}(x)^{-1}x^{C_{0}}\exp(-\Delta_0/x)$
 is a solution of the system $\disp\frac{d}{dx}X=G(x)X$. In the particular case where $G=A_\phi$, the
first
 line of $Y_{0}(x)^{-1}x^{C_{0}}\exp(-\Delta_0/x)$ form a basis of solutions of
$\phi$ at $0$.

$\Delta_0=0$ means that $NR(\phi)$ has non-positive slopes. In this
case, $\mathcal{M}_\phi$ and $\phi$ are both called regular at $0$,
the solution of $\phi$ at $0$ are called logarithmic, and one
verifies that the eigenvalues of $C_0$ modulo ${\bZ}$ depend only on
$\mathcal{M}_\phi$ and are called \emph{exponents} of
$\mathcal{M}_\phi$ and $\phi$ at $0$. According to what precedes, if
the  slopes of $NR(\phi)$ at infinity are in $\{0,1\}$, there exist
a finite extension $K'$
 of $K$, a matrix $Y_{\infty}(x)\in {\GL}_\mu(K'((x)))$, called a
reduction matrix of $A_{\phi}$ (or of $\phi$)
  at infinity, an upper triangular
matrix $C_{\infty}\in
  {\M}_{\mu}(K')$ and a diagonal matrix  $\Delta_\infty$ commuting
with $C_\infty$ such that
$Y_\infty(\frac{1}{x})[A_{\phi}(x)]=-\Delta_\infty-\frac{1}{x}C_{\infty}$.
In this case,
$Y_\infty(\frac{1}{x})^{-1}(\frac{1}{x})^{C_\infty}\exp(-\Delta_\infty
x)$ is a solution of the system $\disp\frac{d}{dx}X=A_{\phi}(x)X$ at
infinity. The exponents of ${\mathcal M}_\phi$ at infinity are those
of ${\mathcal M}_{\phi_\infty}$ at $0$.

By extension, we attribute to $\phi$ the properties that
 ${\mathcal M}_{\phi}$ has. Then one observes, from \S 2.1, that: \\
 $\phi$ is regular at $0$ (resp. infinity) if and only if $\phi^{*}$ is regular at
 the same point,
in which case, the exponents of $\phi^{*}$ at $0$ (resp. at
infinity) are those of $\phi$  but with the opposite sign.

In the sequel, we assume $K$ is sufficiently large so that we can
take $K'=K$. Also, for any matrix $Y$ of ${\M}_\mu(K((x)))$ and any
finite place  $v$ of $V_{0}$, we denote by $r_{v}(Y)$ the upper
bound of the reals $r>0$ for which all entries of $Y$ are analytic
in the punctured open disc $D(0,r^{-})\setminus\{0\}$ of $K_{v}$. If
$Y\in {\GL}_\mu( K((x)))$, we put $R_v(Y)=\min
(r_v(Y),r_v(Y^{-1}))$.
 We end this paragraph with the following
 result due to F. Baldassarri (See Theorem 2 of \cite[III]{Ba}):
\begin{proposition}
 If $Y(x)$ is a reduction matrix of a $K(x)$-module $\mathcal M$
 at $0$ or at infinity, then $R_{v}(Y)$ is non-zero for each finite place $v$
of $V_0$.
\end{proposition}
\subsection{2.4 $\mathcal{E}$-functions}\ \\
A formal power series $f=\sum_{n\ge 0}a_{n}x^{n}\in K[[x]]$ is said
to be  an $\mathcal{E}$\emph{-function}, if the power series
$\sum_{n\ge 0}\frac{a_{n}}{n!}x^{n}$ is a $G$-function.

This definition is motivated by the fact that all power series
occurring in the solutions of the $E$-operators at infinity are
$\mathcal{E}$-functions (see Theorem 2.3 below). A simple example of
these power series is the Euler series : $\sum_{n\ge
0}(-1)^{n}n!x^n$.

We suppose in the sequel that $K$ contains all the coefficients of
 $G$-functions,   $E$-functions and
$\mathcal{E}$-functions that we shall meet thereafter.

The  \emph{Pochhammer symbol}   $(\alpha)_n$ stands for $(\alpha)_n
= \alpha(\alpha+1)\cdots (\alpha+n-1).$ With this notation, the
theorem 2 of \cite{Cl} shows that for any finite place $v$ of $K$
above a prime number $p$, if $\alpha$ is either an integer $\ge 1$
or a non-integer rational number of denominator prime to $p$, then
\begin{equation}
 \lim_{n\to \infty} |(\alpha)_n|_v^{1/n}=\pi_v =p^{-1/(p-1)}.
\end{equation}
For the special case $\alpha=1$, we get
\begin{equation}
\disp\lim_{n\to\infty} |n!|_n^{1/n}=\pi_v
\end{equation}
Combining this equality with the remark below, we find that any
$\mathcal{E}$-function $f$ satisfies
\begin{equation}
\prod_{v\in V_{0}}\min(r_{v}(f)\pi_{v},1)\ne 0
\end{equation}
\begin{remark}{\rm \cite[p 126]{A1}
 If $g$ is a $G$-function, then  $\prod_{v\in
V_{0}}\min(r_{v}(g),1)\ne 0$.}
\end{remark}
\subsection{2.5 $G$-operators and $E$-operators}\ \\
We will give here an equivalent definition of the $G$-operators,
called the \emph{local Bombieri
property} \cite{A1}, which will be useful for the proof of our main theorem.\\
\textbf{Definitions.}  1) An operator $\phi$ of $K[x,d/dx]$ of rank
$\nu$ is said to be a $G$-operator if the differential system
$dX/dx=A_\phi X$ (where $A_\phi$ is the companion matrix of $\phi$
defined in \S2.1) has a solution at $0$ of the form $Y(x)x^C$ where
$Y(x)$ is a $\nu\times\nu$ invertible matrix with entries in
$K((x))$ such that $\prod_{v\in V_{0}}\min(R_{v}(Y_{\phi}),1)\ne 0$,
and where $C$ is a $\nu\times\nu$  matrix with entries in $K$ and
with eigenvalues in ${\bQ}$.\\
2) An operator $\psi\in K[x,d/dx]$ is said to be $E$-operator if it
is the Fourier transform of a certain $G$-operator.

Combining the condition of Bombieri, mentioned in the introduction,
with the properties of the function of generic radius of
convergence, we obtain that if $\phi\in K[x,d/dx]$ is a $G$-operator
then:\\
1) $\overline{\phi}$ and $\phi^*$ are also $G$-operators.\\
2) $\phi$ has only regular singularities with  rational exponents.

From the fact that $
\mathcal{F}(\phi^*)=\Big(\overline{\mathcal{F}(\phi)}\Big)^*$ (cf.
\cite[V.3.6]{Ma}), the first statement implies that, if $\psi$ is an
$E$-operator then $\overline{\psi}$ and  $\psi^*$ are also
$E$-operators. The second statement means that  the Newton-Ramis of
any $G$-operator has non slope other than $0$ and $\infty$, and
hence, from \S2.2, that the slopes of Newton-Ramis of any
$E$-operator are in $ \{0,-1\}$.

The following theorem, due to Andr\'e, describes the nature of
solutions of the $E$-operators at $0$ and at infinity:
\begin{theorem} \cite{A2} Let $\psi$ be an $E$-operator of rank $\mu$,
then :
\begin{itemize}
\item[(1)] The slopes of $NR(\psi)$ lie in  $\{-1,0\}$; \item[(2)]
$\psi$ admits a basis of solutions at $0$ of the form
$$(F_{1},...,F_{\mu})x^{\Gamma_0}$$
 where the $F_{i}$ are
$E$-functions, where $\Gamma_0$ is a $\mu\times\mu$ upper triangular
 matrix  with elements in $\mathbb{Q}$;
\item[(3)] $\psi$ admits a basis of solutions at infinity of the form
$$\Big(f_{1}\Big(\frac{1}{x}\Big),...,f_{\mu}\Big(\frac{1}{x}\Big)\Big)\Big(\frac{1}{x}\Big)^{\Gamma}\exp(-\Delta
x)$$
 where the $f_{i}$ are
$\mathcal{E}$-functions, where $\Gamma$ is a $\mu\times\mu$ upper
triangular  matrix  with elements in $\mathbb{Q}$, and where
$\Delta$ is a $\mu\times\mu$ diagonal matrix with elements in $K$
which commutes with $\Gamma$.
\end{itemize}
\end{theorem}
\section{The main theorem}
Before stating the main theorem, we recall that for a given
differential operator $\psi$, $A_\psi$ denotes its companion matrix
(see \S 2.1).
\begin{theorem}
Let $\psi$ be a differential operator of $ K[x,d/dx]$. Then $\psi$
is an $E$-operator if and only if $\psi$ satisfies the following
conditions :
\begin{itemize}
\item[(1)]
the coefficients of $\psi$ are not all in K;\item[(2)]
 the slopes of
$NR(\psi)$ lie in $\{-1,0\}$;
\item[(3)] the differential system $dZ/dx=A_{\psi}Z$  has a solution
of the from
$$Y(\frac{1}{x})(\frac{1}{x})^\Gamma\exp(-\Delta x),$$
where $Y(x)$ is a $\mu\times\mu$ invertible matrix  with entries in
$K((x))$ such that $\prod_{v\in V_{0}} \min( R_{v}(Y)\pi_{v},1)\ne
0,$ where $\Gamma$ is a $\mu\times\mu$  matrix with entries in $K$
and with eigenvalues in $\mathbb{Q}$, and where $\Delta$ is a
$\mu\times\mu$ diagonal matrix with entries in $K$ which commutes
with $\Gamma$.
\end{itemize}
\end{theorem}
The first condition of this theorem is necessary by definition of
the $E$-operators. Theorem 2.3 above shows that the seconde one is
also necessary. In \S3.2 we will prove that the third one is also
necessary. The fact that these conditions are sufficient is
postponed to section 6. In the following paragraph, we give some
preliminary results which will be useful in the rest of this paper.
\subsection{3.1. Preliminary results}\ \\
Throughout this paragraph, $\phi=a_{\mu}(d/dx)^\mu+\ldots+a_{0}$
 denotes a differential operator of $ K[x,d/dx]$ of rank $\mu\in
{\bZ}_{>0}$, $\overline{\phi}$ denotes  the differential operator
obtained from $\phi$ by change of variable $x\to -x$, $\Gamma_1$ and
$\Gamma_2$ denote two $\mu\times\mu$ matrices
 with entries in $K$, $\Delta_1$  and $\Delta_2$ denote two $\mu\times\mu$ diagonal
matrices  with entries in $K$ such that
$\Gamma_1\Delta_1=\Delta_1\Gamma_1$ and
$\Gamma_2\Delta_2=\Delta_2\Gamma_2$, and
$y_1,\ldots,y_\mu,z_1,\ldots,z_\mu$  denote power series of
$K((x))$.
\begin{lemma} Let $G$ be a $\mu\times\mu$ matrix with entries in $K$, and let $Y_1$ and $Y_2$ be
two matrices of ${\GL}_\mu(K((x)))$ such that
$Y_1[G]=\disp\frac{\Delta_1}{x^2}+\frac{\Gamma_1}{x}$  and
$Y_2[G]=\disp\frac{\Delta_2}{x^2}+\frac{\Gamma_2}{x}$. Then,
\begin{itemize}
\item[(1)] the matrices $\Delta_1$ and $\Delta_2$ are similar; \item[(2)]
$Y_1Y_2^{-1}[\frac{\Gamma_2}{x}]=\frac{\Gamma_1}{x}$ and
$Y_1Y_2^{-1}\in {\GL}_\mu(K[x,1/x])$;
\item[(3)] the eigenvalues of $\Gamma_1$ coincide, modulo ${\bZ}$,
with those of $\Gamma_2$.
\end{itemize}
\end{lemma}
\begin{proof} Let $a=(a_1,\ldots,a_\mu)\in K^\mu$. Put, for
$i=1,2$,
$$\mathrm{E}_i(a)=\{\mathrm{v}\in K^\mu\;|\;
\Delta_i\mathrm{v}=a_j\mathrm{v},\;\; 1\le j\le \mu\},$$ and
$$\Sigma_{i}=\{a\in K^\mu\;|\; \mathrm{E}_i(a)\ne 0\}.$$ Then
$$K^\mu=\bigoplus_{a\in \Sigma_{i}}\mathrm{E}_i(a).$$  Moreover, $\Gamma_i$ commutes
with the  projection $K^\mu\longrightarrow \mathrm{E}_i(a)$. Thus,
$\Gamma_i$ can be written as $$ \Gamma_i=\bigoplus_{a\in
\Sigma_{i}}\Gamma_i(a),\;\;\;\text{where}\;\;\;\;\; \Gamma_i(a)\in
{\M}_{\dim_K(\mathrm{E}_i(a))}(K).$$ In addition, by
 hypothesis, we have
$$Y_1Y_2^{-1}[\disp\frac{\Delta_2}{x^2}+\frac{\Gamma_2}{x}]=\disp\frac{\Delta_1}{x^2}+\frac{\Gamma_1}{x}.$$
According to  Proposition 6.4 of \cite{BV}, we find \\
1) The matrices $\Delta_1$ and $\Delta_2$ are similar,\\
2)
$\Sigma:=\Sigma_{1}=\Sigma_{2},\;\;\;\dim_K(\mathrm{E}_1(a))=\dim_K(\mathrm{E}_2(a))\;\;\;\text{for
any}\;\;a\in \Sigma,$\\
3) $Y:=Y_1Y_2^{-1}=\bigoplus_{a\in \sum}Y(a),$ such that $Y(a)\in
{\M}_{\dim_K(\mathrm{E}_1(a))}(K((x)))$ and
$Y(a)\Big[\Gamma_2(a)/x\Big]=\Gamma_1(a)/x$ for any $a\in
\Sigma$.\\

Thus, $Y\Big[\Gamma_2/x\Big]=\Gamma_1/x$ and hence the eigenvalues
of $\Gamma_1$ coincide, modulo ${\bZ}$, with those of $\Gamma_2$
(cf. \cite[III.8]{DGS}). Moreover, for any $a\in \Sigma$, we have
$$x\disp\frac{d}{dx}Y(a)=\Gamma_1(a)Y(a)-Y(a)\Gamma_2(a).$$ Therefore, if we write
$Y(a)= \disp\sum_{m\in\mathbb{Z}}Y(a)_{m}x^{m} $, we obtain for any
$m\in \mathbb{Z}\setminus\{0\}$,
$$mY(a)_{m}=\Gamma_1(a)Y(a)_{m}-Y(a)_{m}\Gamma_2(a).$$
But the eigenvalues of the maps
\begin{eqnarray*}
   T_{m}(a)\;\colon {\M}_\mu(\overline{{\bQ}}) &\longrightarrow &{\M}_\mu(\overline{{\bQ}}) \\
           X &\longmapsto &\Gamma_1(a)X-X\Gamma_2(a)-m X
\end{eqnarray*}
are of the form $\lambda(a)-\gamma(a)-m$ where $\lambda(a)$ and
$\gamma(a)$
 are  respectively  eigenvalues of $\Gamma_1(a)$ and of $\Gamma_2(a)$.
This means that $T_m(a)$ is invertible, except perhaps, for a finite
set of integers $m$. Hence, $Y(a)_{m}$ is zero except
 for a finite set of integers $m$ and the conclusion  follows.
\end{proof}
\begin{corollary}
Let $(y_1,\ldots,y_\mu)x^{\Gamma_1}\exp(\Delta_1/x)$  be a  basis of
solutions of $\phi$ at $0$. Then, $\phi$ has a basis of solutions at
$0$ of the form
$(\widetilde{y}_1,\ldots,\widetilde{y}_\mu)x^{\widetilde{\Gamma}}\exp(\widetilde{\Delta}/x)=(\xi_1,\ldots,\xi_{\mu})$
where:
\begin{itemize}
\item[(1)]$\widetilde{y}_1,\ldots,\widetilde{y}_\mu$ are formal power series of $K((x))$
 and are $K[x,1/x]$-linear combinations of
$y_1,\ldots,y_\mu$ and of their derivatives; \item[(2)]
$\widetilde{\Gamma}$ is a $\mu\times\mu$ matrix, in Jordan form,
whose entries lie in $K$ and whose eigenvalues coincide, modulo
${\bZ}$, with those of $\Gamma_1$;
\item[(3)] $\widetilde{\Delta}$ is a $\mu\times\mu$ diagonal matrix similar to $\Delta_1$.
\end{itemize}
Moreover, if $\gamma_1,\ldots,\gamma_\mu$ denote the eigenvalues of
$\Gamma_1$ and $\delta_1,\ldots,\delta_\mu$ denote the diagonal
terms of $\Delta_1$, then $\xi_1,\ldots,\xi_{\mu}$ lie in
$$\Big<\widetilde{y}_i\;x^{\gamma_j}\;(\ln
x)^{k-1}\exp(\delta_\ell/x),\;\; 1\le i,j,k,\ell\le
\mu\Big>_{K[x,1/x]}$$
\end{corollary}
\begin{proof} Let  $W_1$ be the wronskian matrix
of  $(y_1,\ldots,y_\mu)x^{\Gamma_1}\exp(\Delta_1/x)$. Thus, $W_1$ is
a solution of the system $dX/dx=A_{\phi}X$. Moreover, $W_1$ can be
written  of the form $Y_1x^{\Gamma_1}\exp(\Delta_1/x)$, where $Y_1$
is a matrix of ${\GL}_\mu(K((x)))$ whose entries are
$K[x,1/x]$-linear combinations of $y_1,\ldots,y_\mu$ and of their
derivatives. Therefore, $Y_1^{-1}[A_\phi]=-\Delta_1/x^2+\Gamma_1/x$.
The Turrittin-Levelt decomposition states, in this case, that there
exists  a $\mu\times\mu$ invertible matrix
$\widetilde{Y}=(\widetilde{y}_{ij})\in {\GL}_{\mu}(K((x)))$,  a
$\mu\times\mu$ matrix $\widetilde{\Gamma}$ in Jordan form with
entries in $K$, and a $\mu\times\mu$ diagonal matrix
$\widetilde{\Delta}=(\widetilde{\delta}_{ij})$ with entries in $K$
commuting with $\widetilde{\Gamma}$ such that
$\widetilde{Y}^{-1}[A_{\phi}]=-\widetilde{\Delta}/x^2+\widetilde{\Gamma}/x$.
Hence, by previous Lemma, the matrices $\Delta_1$ and
$\widetilde{\Delta}$ are similar, the eigenvalues of
$\widetilde{\Gamma}$ coincide, modulo ${\bZ}$, with those of
$\Gamma_1$, and there exists $L\in {\GL}_{\mu}(K[x,1/x])$ such that
$\widetilde{Y}=LY_1$. In particular, the entries
$(\widetilde{y}_{ij})$ of $\widetilde{Y}$ are $K[x,1/x]$-linear
combinations of $y_1,\ldots,y_\mu$ and of their derivatives. In
addition, since the matrix
$\widetilde{Y}x^{\widetilde{\Gamma}}\exp(\widetilde{\Delta}/x)$ is a
solution of the system $dX/dx=A_{\phi}X$, it is the wronskian matrix
of the $\mu$-tuple
$(\widetilde{y}_{11},\ldots,\widetilde{y}_{1\mu})x^{\widetilde{\Gamma}}\exp(\widetilde{\Delta}/x)$.
Thus, the coefficients of
$(\widetilde{y}_{11},\ldots,\widetilde{y}_{1\mu})x^{\widetilde{\Gamma}}\exp(\widetilde{\Delta}/x)$
form a basis of solutions of $\phi$ at $0$. Hence,  by putting
$\widetilde{y}_{i}=\widetilde{y}_{1i}$ for $i=1,\ldots,\mu$, we find
that
$(\widetilde{y}_{1},\ldots,\widetilde{y}_{\mu})x^{\widetilde{\Gamma}}\exp(\widetilde{\Delta}/x)$
is a basis of solutions of $\phi$ at $0$ which meets the conditions
(1), (2) and (3) of Corollary 3.3. On the other hand, by hypotheses,
$\widetilde{\Gamma}=(\widetilde{\gamma}_{ij})$ is of the form $D+N$
where $D$ is a diagonal matrix and $N$ is a nilpotent upper
triangular matrix such that $DN=ND$ and $N^\mu=0$. Thus, $$\disp
x^{\widetilde{\Gamma}}= x^{D+N}=x^{D}\sum_{ 0\le k\le
\mu-1}\frac{N^{k}}{k!}(\ln x)^k =x^{D}+x^{D}\sum_{1\le k\le
\mu-1}\disp\frac{N^{k}}{k!}(\ln x)^k.$$ Therefore, $
\xi_{1}=\widetilde{y}_{1}x^{\widetilde{\gamma}_{11}}\exp(\widetilde{\delta}_{11}/x),
$ and for all $2\le i\le \mu$, $$
\xi_{i}=\Big(\widetilde{y}_{i}x^{\widetilde{\gamma}_{ii}}+\disp\sum_{j=1}^{
i-1}\widetilde{y}_{j}x^{\widetilde{\gamma}_{jj}}\sum_{1\le k\le
\mu-1}\disp\frac{(N^{k})_{ji}}{k!}(\ln
x)^k\Big)\exp(\widetilde{\delta}_{ii}/x),
$$ since $(N^{k})_{ji}=0$  for all  $k\ge 1$ and all
$1\le i\le j\le \mu$. Hence, the last statement of the corollary
results from the fact that $\Delta_1$ and $\widetilde{\Delta}$ are
similar and that $\widetilde{\Gamma}$ and $\Gamma_1$ have the same
eigenvalues modulo ${\bZ}$.
\end{proof}
\begin{corollary}
Let $(y_1,\ldots,y_\mu)x^{\Gamma_1}\exp(\Delta_1/x)$ and
$(z_1,\ldots,z_\mu)x^{\Gamma_2}\exp(\Delta_2/x)$ be respectively
bases of solutions of $\phi$ and $\phi^*$ at $0$. Then, the
differential system $dX/dx=A_{\phi}X$ has a solution of the form
$Y(x)x^{\Gamma_1}\exp(\Delta_1/x)$, where $Y$ is  a $\mu\times\mu$
invertible matrix such that the entries of
 $Y$ $($resp. of $Y^{-1})$ are $K[x,1/x]$-linear
 combinations of $y_1,\ldots,y_\mu\;($resp.
 of $z_1,\ldots,z_\mu)$ and of their derivatives. Moreover, the matrices
 $\Delta_1$ and $-\Delta_2$ are similar, and
 the eigenvalues of $\Gamma_1$ are those of $-\Gamma_2$ modulo ${\bZ}$.
 \end{corollary}
 \begin{proof} Let us write $\phi=a_{\mu}(d/dx)^\mu+\ldots+a_{0}$. Since
$(y_1,\ldots,y_\mu)x^{\Gamma_1}\exp(\Delta_1/x)$  is a basis of
solutions of $\phi$ at $0$, the wronskian matrix $W$ of the
$\mu$-tuple $(y_1,\ldots,y_\mu)x^{\Gamma_1}\exp(\Delta_1/x)$  is
then a solution of the system $dX/dx=A_{\phi}X$. Moreover, $W$ can
be written in the form $Yx^{\Gamma_1}\exp(\Delta_1/x)$ where $Y$ is
a matrix of ${\GL}_\mu(K((x)))$ whose entries are $K[x,1/x]$-linear
combinations of $y_1,\ldots,y_\mu$ and of their derivatives. Then we
have, $Y^{-1}[A_\phi]=-\Delta_1/x^2+\Gamma_1/x$, which means that
$Y^{-1}$ is a reduction matrix of $\phi$ at $0$, or also
\begin{eqnarray}
^{T}Y[-^{T}A_{\phi}]&=&^{T}\Delta_1\disp\frac{1}{x^2}-^{T}\Gamma_1\disp\frac{1}{x}.
\end{eqnarray}
In addition, the $\mu$-tuple $a_{\mu}(x)
(z_1,\ldots,z_\mu)x^{\Gamma_2}\exp(\Delta_2/x)$  is a  basis of
solutions of $\phi^{*}a_{\mu}^{-1}=(a_{\mu}^{-1}\phi)^{*}$
 at $0$.
  Therefore, the matrix $U$  whose rows $u_1,\ldots,u_\mu$  are
 defined recursively by
\begin{eqnarray*}
u_\mu&=&a_{\mu}(x)
(z_1,\ldots,z_\mu)x^{\Gamma_2}\exp(\Delta_2/x)\\
u_{\mu-i}&=&\disp\frac{a_{\mu-i}(x)}{a_{\mu}(x)}u_{\mu}-\disp\frac{d}{dx}u_{\mu-i+1}\;\;\;(1\le
i\le \mu-1),
\end{eqnarray*}
 is a solution of the system $dX/dx=-^{T}A_{\phi}X.$ Moreover, $U$  may be
written of the form $Zx^{\Gamma_2}\exp(\Delta_2/x)$, where $Z$  is
an invertible matrix  $\mu\times \mu$ whose entries are
$K[x,1/x]$-linear combinations of $z_1,\ldots,z_\mu$ and of their
derivatives. Thus, we have
\begin{eqnarray}
Z^{-1}[-^{T}A_{\phi}]&=&-\disp\frac{1}{x^2}\Delta_2+\disp\frac{1}{x}\Gamma_2.
\end{eqnarray}
Thus, by formulae (3.1), (3.2) and Lemma 3.2, the matrices
$^T\Delta_1(=\Delta_1)$ and $-\Delta_2$ are similar, the eigenvalues
of $^T\Gamma_1$(which are also those of $\Gamma_1$) are those of
$-\Gamma_2$ modulo ${\bZ}$ and there exists $L\in
{\GL}_\mu(K[x,1/x])$ such that $^{T}Y=LZ^{-1}$. Consequently, the
entries of $Y^{-1}$ are $K[x,1/x]$-linear combinations of
$z_1,\ldots,z_\mu$ and of their derivatives. The conclusion follows.
\end{proof}
\begin{lemma}
Let $(y_{1}(x),\ldots,y_{\mu}(x))x^{\Gamma_1}\exp(\Delta_1/x)$ be a
basis of solutions of $\phi$ at $0$. Then
$(y_{1}(-x),\ldots,y_{\mu}(-x))x^{\Gamma_1}\exp(-\Delta_1/x)$ is a
basis of solutions of $\overline{\phi}$ at $0$.
\end{lemma}
\begin{proof}
Let $W$ be the wronskian matrix of
$(y_{1}(x),\ldots,y_{\mu}(x))x^{\Gamma_1}\exp(\Delta_1/x)$. Thus $W$
can be written in the form $Y(x)x^{\Gamma_1}\exp(\Delta_1/x)$, where
$Y(x)$ is a $\mu\times\mu$ invertible matrix with entries in
$K((x))$. Thus $Y^{-1}(x)[A_\phi]=-\Delta_1/x^2+\Gamma_1/x$. By
change of variable $x\to -x$, we find
$$Y^{-1}(-x)A_{\phi}(-x)Y(-x)+Y^{-1}(-x)\frac{d}{dx}(Y(-x))=-\frac{\Delta_1}{x^2}-\frac{\Gamma_1}{x}.$$
Thus $Y(-x)x^{\Gamma_1}\exp(-\Delta_1/x)$ is solution of the system
$\frac{d}{dx}X=-A_{\phi}(-x)X$. Consequently, the $\mu$-tuple
$(y_{1}(-x),\ldots,y_{\mu}(-x))x^{\Gamma_1}\exp(-\Delta_1/x)$ is a
basis of solutions of $\overline{\phi}$ at $0$.
\end{proof}
\begin{lemma}Let $\phi=a_{\mu}(d/dx)^\mu+\ldots+a_{0}\in K[x,d/dx]$.
 Let $W$  be a $\mu\times\mu$ invertible matrix with entries in some
Picard-Vessiot extension of $K$. If $W$ is a  solution of
$dX/dx=A_\phi X$ at $0$. Then the elements of $\mu-$th row of
$a_{\mu}^{-1}(^TW^{-1})$ form a basis of solutions of $\phi^*$ at
$0$.
\end{lemma}
\begin{proof}
Let  $W_1,\ldots,W_\mu$ denote the rows of $^TW^{-1}$. Since $W$ is
a  solution of $dX/dx=A_\phi X$ at $0$, these rows are then related
by
$$W_{\mu-i}=\frac{a_{\mu-i}}{a_{\mu}}W_{\mu}-\frac{d}{dx}W_{\mu-i+1}, \;\;\; (1\le i\le \mu-1),$$
and the elements of  $W_\mu$ are solutions of
$\phi^*a_\mu^{-1}=(a_\mu^{-1}\phi)^*$ at $0$. We get therefore, by
 induction on the index  $i$,
\begin{equation*}
\label{adj}W_i\in
\Big<W_\mu,\ldots,\Big(\frac{d}{dx}\Big)^{\mu}(W_\mu)\Big>_{K[x,a_{\mu}^{-1}]},\;\;\;
(1\le i\le \mu-1).
\end{equation*}
In addition, since $W$ is an invertible matrix, the elements of
$W_\mu$ are then linearly independent over $K$, and hence they form
a basis of solutions of $\phi^*a_\mu^{-1}$ at $0$, and the
conclusion follows.
\end{proof}
\begin{lemma}  Let $Y(x)x^{\Gamma_1}\exp(\Delta_1/x)$ be a solution
of $dX/dx=A_{\phi}(x)X$ at $0$, where $Y(x)=(y_{ij}(x))\in
{\GL}_{\mu}(K((x)))$. Assume that the leading coefficient $a_{\mu}$
of $\phi$ is a monomial. Then there exists $\widetilde{Y}(x)\in
{\GL}_{\mu}(K((x)))$ such that
$\widetilde{Y}(x)x^{-^T\Gamma_1}\exp(\Delta_1/x)$ is a solution
 of $dX/dx=A_{\overline{\phi^*}}(x)X$ at $0$, and
$R_v(\widetilde{Y})\ge R_v(Y)$ for all $v\in V_0$.
\end{lemma}
\begin{proof}  First, it is easy to check that
$Y(x)x^{\Gamma_1}\exp(\Delta_1/x)$ is the wronskian matrix of the
elements  of
$(y_{11}(x),\ldots,y_{1\mu}(x))x^{\Gamma_1}\exp(\Delta_1/x)$ and
that the elements of
$(y_{11}(x),\ldots,y_{1\mu}(x))x^{\Gamma_1}\exp(\Delta_1/x)$ is a
basis of $\phi$ at $0$.\\  In addition,  if we write
$Y^{-1}(x)=(\widetilde{y}_{ij}(x))$, we find, by Lemma 3.6, that
$a_\mu^{-1}(x)((\widetilde{y}_{1\mu}(x),\ldots,\widetilde{y}_{\mu\mu}(x))x^{-^T\Gamma_1}\exp(-\Delta_1/x)$
is a basis of solutions $\phi^*$ at $0$. According to Lemma 3.5, the
$\mu$-tuples
$(y_{11}(-x),\ldots,y_{1\mu}(-x))x^{\Gamma_1}\exp(-\Delta_1/x)$ and
$a_\mu^{-1}(-x)((\widetilde{y}_{1\mu}(-x),\ldots,\widetilde{y}_{\mu\mu}(-x))x^{-^T\Gamma_1}\exp(\Delta_1/x)$
are respectively bases of solutions of
$\overline{\phi}=(\overline{\phi^*})^*$ and $\overline{\phi^*}$ at
$0$. Finally, since $a_\mu^{-1}$ is a monomial, Lemma 3.4 states
that $dX/dx=A_{\overline{\phi^*}}(x)X$ has a solution at $0$ in the
form $\widetilde{Y}(x)x^{-^T\Gamma_1}\exp(\Delta_1/x)$, where
$\widetilde{Y}(x)$ is a $\mu\times\mu$ invertible matrix such that
the entries of
 $\widetilde{Y}(x)$ $($resp. of $\widetilde{Y}^{-1}(x))$ are $K[x,1/x]$-linear
 combinations of
$\widetilde{y}_{1\mu}(-x),\ldots,\widetilde{y}_{\mu\mu}(-x)$ (resp.
 of $y_{11}(-x),\ldots,y_{1\mu}(-x)$) and of their derivatives.
 Hence, for all $v\in V_0$, we have
 $R_{v}(\widetilde{Y})=\min_{1\le i,j\le\mu}\{r_{v}(y_{1i}),r_{v}(\widetilde{y}_{j\mu})\}\ge
 R_{v}(Y)$. The conclusion follows.
 \end{proof}
\subsection{3.2. Necessary conditions}\ \\
We conclude this section by proving that the seconde condition of
Theorem 3.1 is necessary:
\begin{theorem}Let $\psi$ be an
$E$-operator of $K[x,d/dx]$ of rank $\mu$. Then, the differential
system $d/dxZ=A_{\psi}Z$ has a solution of the from
$$Y(\frac{1}{x})(\frac{1}{x})^\Gamma\exp(-\Delta x),$$
where $Y(x)$ is an $\mu\times\mu$ invertible matrix  with entries in
$K((x))$ such that $\prod_{v\in V_{0}} \min( R_{v}(Y)\pi_{v},1)\ne
0,$ where $\Gamma$ is a $\mu\times\mu$ upper triangular  matrix with
entries in $\mathbb{Q}$, and where $\Delta$ is a $\mu\times\mu$
diagonal  matrix with entries in $K$ which commutes with $\Gamma$.
\end{theorem}
\begin{proof} According to \S2.5, the operator $\psi^*$ is also an $E$-operator.
Combining this, with Theorem 3.1 and Corollary 3.4 (applied at
infinity), we observe that the differential system $d/dxZ=A_{\psi}Z$
has a solution of the from
$$Y(\frac{1}{x})(\frac{1}{x})^\Gamma\exp(-\Delta x),$$
where $Y(x)$ is a $\mu\times\mu$ invertible matrix such that the
entries  of $Y(x)=(y_{ij})$ and those of
$Y(x)^{-1}=(\widetilde{y}_{kl})$ are $K[x,1/x]$-linear
 combinations of $\mathcal{E}$-functions and of their derivatives,
 where $\Gamma$ is a $\mu\times\mu$ upper triangular  matrix
with entries in $\mathbb{Q}$, and where $\Delta$ is a $\mu\times\mu$
diagonal matrix with entries in $K$ which commutes with $\Gamma$.
Thus, by (2.4), we have
\begin{eqnarray*}
\prod_{v\in V_{0}}\min(R_{v}(Y_\psi)\pi_{v},1)
&=&\prod_{v\in V_{0}}\Big(\min(\min_{i,j}(r_{v}(y_{ij})\pi_{v}),\min_{k,l}(r_{v}(\widetilde{y}_{kl})\pi_{v}),1)\Big)\\
                 &\ge &\prod_{v\in
V_{0}}\Big(\prod_{ij}\min(r_{v}(y_{ij})\pi_{v},1)\prod_{kl}\min(r_{v}(\widetilde{y}_{kl})\pi_{v},1)\Big)\\
&\ge& \prod_{ij}\prod_{v\in
V_{0}}\min(r_{v}(y_{ij})\pi_{v},1)\\
&&\;\;\;\times\;\prod_{kl}\prod_{v\in
V_{0}}\min(r_{v}(\widetilde{y}_{kl})\pi_{v},1)\ne 0.
\end{eqnarray*}
\end{proof}
In the sequel, we fix an embedding of $K$ into ${\bC}$.
\section{The Laplace transform}
\subsection{4.1 The Laplace transform $\mathcal L$}\ \\
In this paragraph, we summarize main properties of the formal
Laplace transform due in part to Y. Andr\'e.

Let $\alpha$ be an element of $K$ with real part  $>-1$, let $k$ and
$n$ be two nonnegative integers, and let $h_{\alpha,k}$ denote the
function defined by $h_{\alpha,k}(x)=x^{\alpha}(\ln x)^{k}$; $x>0$.
The standard Laplace transform of $h_{\alpha,0}$, denoted ${\mathcal
L}(h_{\alpha,0})$, is given by (cf. \cite[2.30]{DP})
\begin{eqnarray}
{\mathcal L}(h_{\alpha,0})(z)&=& \int_0^\infty e^{-zx}x^\alpha
dx=\Gamma(\alpha+1)z^{-\alpha-1}.
\end{eqnarray}
This implies in particular,
$$(\frac{d}{d\alpha})^{k}\Big(\Gamma(\alpha+1)z^{-\alpha-1}\Big)=
\int_{0}^{\infty}e^{-zx}x^{\alpha} (\ln x)^{k}dx={\mathcal
L}(h_{\alpha,k})(z).$$ Leibniz formula gives
$$(\frac{d}{d\alpha})^{k}\Big(\Gamma(\alpha+1)z^{-\alpha-1}\Big)=
\sum_{j=0}^{k}\binom{k}{j}\Gamma^{(j)}(\alpha+1)z^{-\alpha-1}(-1)^{k-j}(\ln
z)^{k-j}.$$ Thus,
\begin{eqnarray}
{\mathcal L}(h_{\alpha,k})(z)=
\sum_{j=0}^{k}\binom{k}{j}\Gamma^{(j)}(\alpha+1)z^{-\alpha-1}(-1)^{k-j}(\ln
z)^{k-j}.
\end{eqnarray}
From the fact that $\Gamma(\alpha+1)=\alpha\Gamma(\alpha),$ we
obtain, by induction on $j\ge 1$, the  following relations
$$\Gamma'(\alpha+1)=\Gamma(\alpha)+\alpha\Gamma'(\alpha)\;\text{and}\;\;\;
\Gamma^{(j+1)}(\alpha+1)=j\Gamma^{(j)}(\alpha)+\alpha\Gamma^{(j+1)}(\alpha),$$
which implies
\begin{eqnarray}
{\mathcal L}(h_{\alpha,k})(z)&\in&
z^{-\alpha-1}<\Gamma(\alpha),\ldots,\Gamma^{(k)}(\alpha)>_{{\bQ}[\alpha,\ln
z]},
\end{eqnarray}
and, in the case where $\alpha$ is a non-zero positive integer,
gives
\begin{eqnarray}
\Gamma'(\alpha)&=&\alpha!.
\end{eqnarray}
On the other hand, the function $h_{\alpha,k}$ satisfies the
following equalities (cf. \cite[2.21, 2.40]{DP}),
\begin{equation}
\begin{aligned}
\frac{d}{dz}{\mathcal L}(h_{\alpha,k})(z)&= {\mathcal
L}(-xh_{\alpha,k})(z).
\end{aligned}
\end{equation}
\begin{equation}
\begin{aligned}{\mathcal L}(\frac{d}{dx}h_{\alpha,k})(z)&= z {\mathcal
L}(h_{\alpha,k})(z)+\lim_{x\longrightarrow
0^{+}}h_{\alpha,k}(x)\\
&=z {\mathcal L}(h_{\alpha,k})(z)\;\;\;\text{if}\;\;\; \Re
e(\alpha)>0.
\end{aligned}
\end{equation}
To extend the Laplace transform  ${\mathcal L}$ of $h_{\alpha,k}$ to
any  $\alpha$, we have to introduce the finite parts  of
$h_{\alpha,k}$ in
the following  manner:\\
Putting $$\Phi(x,\alpha,k)=\int h_{\alpha,k}(x)dx,$$ we find
therefore
\begin{eqnarray*}
\Phi(x,\alpha,k)&=&\frac{x^{\alpha+1}}{\alpha+1}\sum_{\ell=0}^{k}\frac{(-1)^{k-\ell}k!}{(\alpha+1)^{k-\ell}\ell!}(\ln
x)^{\ell}\;\;\;\;\;\; \text{if}\;\;\;\; \alpha\ne -1\\
\Phi(x,-1,k)&=&\frac{(\ln x)^{k+1}}{k+1}.
\end{eqnarray*}
The  finite part of the integral
$\disp\int_{0}^{x}\sum_{\alpha,k}\lambda_{\alpha,k}h_{\alpha,k}(t)dt$,
where the $\lambda_{\alpha,k}$ are complex numbers, is defined, for
$x>0$, by
\begin{equation*}
\begin{aligned}
\text{p.f.}\int_{0}^{x}\sum_{\alpha,k}\lambda_{\alpha,k}h_{\alpha,k}(t)dt&=\sum_{\alpha,k}\lambda_{\alpha,k}
\lim_{\epsilon\longrightarrow
0^{+}}\Big(\Phi(\epsilon,\alpha,k)+\int_{\epsilon}^{x}h_{\alpha,k}(t)dt\Big)\\
&=\sum_{\alpha,k}\lambda_{\alpha,k}\Phi(x,\alpha,k).
\end{aligned}
\end{equation*}
With this definition, we get
\begin{equation}
\begin{aligned}
h_{\alpha,k}^{0}:=h_{\alpha,k}-\text{p.f.}\int_{0}^{x}\Big(\frac{d}{dt}h_{\alpha,k}(t)\Big)
dt&= 0\;\;\;\;\;\;\;\; \text{if}\;\;\;\;\;\; (\alpha,k)\ne (0,0),\\
&= h_{0,0}\;\;\;\;\;\;\text{otherwise}.
 \end{aligned}
\end{equation}
Now, fix $\alpha\in{\bC}$, $k\in {\bZ}_{\ge 0}$ and put for $n\in
{\bZ}_{\ge 0}$,
$$F_{n}(x)=\text{p.f.}\int_{0}^{x}\frac{(x-t)^{n}}{n!}h_{\alpha,k}(t)dt.$$
Then,
\begin{equation}
\begin{aligned}
F_{n}(x)&= \sum_{m=0}^{n}
\frac{(-1)^{m}}{m!(n-m)!}\frac{x^{\alpha+n+1}}{m+\alpha+1}
\sum_{\ell=0}^{k}\frac{k!(-1)^{k-\ell}}{\ell!(m+\alpha+1)^{k-\ell}}(\ln
x)^{\ell},\\
    & \hspace{3cm} \text{if}\;\; (-\alpha-1)\ne 0,1,\ldots,n\\
    &=\sum_{m=0\atop m\ne -\alpha-1}^{n}
\frac{(-1)^{m}}{m!(n-m)!}\frac{x^{\alpha+n+1}}{m+\alpha+1}
\sum_{\ell=0}^{k}\frac{k!(-1)^{k-\ell}}{\ell!(m+\alpha+1)^{k-\ell}}(\ln
x)^{\ell}\\
&\hspace{1cm} +
\frac{(-1)^{\alpha+1}}{(-\alpha-1)!(n+\alpha+1)!}\;x^{\alpha+n+1}\;\frac{(\ln
x)^{k+1}}{k+1},\;\;\;\;\text{otherwise}.\\
\end{aligned}
\end{equation}
Moreover, these  functions satisfy, for any integer $n\ge 1$, the
following equality (cf. \cite[5.35]{DP}),
\begin{eqnarray*}
\frac{d}{dx}F_{n}=F_{n-1}.
\end{eqnarray*}
This means, by $(4.6)$, that the function $z^{n+1}{\mathcal
L}(F_{n})$ is  independent of the choice of $n$ for $n\ge -\Re
e(\alpha)-1$. From this remark, we can extend the Laplace transform
to any $\alpha$ by putting
\begin{eqnarray}
{\mathcal L}(h_{\alpha,k})(z)&=& z^{n+1}{\mathcal
L}(F_{n})(z),\;\;\;\text{for}\;\;\; n\ge -\Re e(\alpha)-1.
\end{eqnarray}
For simplicity, we write $\mathcal{L}(.)$ instead of
$\mathcal{L}(.)(z)$. Then using formula $(4.5)$ and linearity of
$\mathcal{L}$, it follows that for any $\alpha$, $k$ and $n\ge -\Re
e(\alpha)-1$,
\begin{equation}
\begin{aligned}
\frac{d}{dz}({\mathcal L}(h_{\alpha,k}))=& (n+1)z^{n}{\mathcal
L}(F_{n})+z^{n+1}\frac{d}{dz}({\mathcal L}(F_{n})\\
 =&(n+1)z^{n}{\mathcal
L}(F_{n})+z^{n+1}{\mathcal L}(-xF_{n})\\
=&(n+1)z^{n}{\mathcal L}(F_{n})-z^{n+1}{\mathcal
L}\Big((n+1)\text{p.f.}\int_{0}^{x}\frac{(x-t)^{n+1}}{(n+1)!}h_{\alpha,k}(t)dt\\
&\hspace{4.2cm}-\text{p.f.}\int_{0}^{x}\frac{(x-t)^{n}}{n!}(-t)h_{\alpha,k}(t)dt\Big)\\
=&(n+1)z^{n}{\mathcal L}(F_{n})-(n+1)z^{n+1}{\mathcal
L}(F_{n+1})+{\mathcal L}(-xh_{\alpha,k})\\
=&{\mathcal L}(-xh_{\alpha,k}).
\end{aligned}
\end{equation}
On the other hand, for $n\ge -\Re e(\alpha)$
\begin{eqnarray*}
{\mathcal L}(\frac{d}{dx}h_{\alpha,k})&=& z^{n+1}{\mathcal
L}\Big(\text{p.f.}\int_{0}^{x}\frac{(x-t)^{n}}{n!}(\frac{d}{dt}h_{\alpha,k}(t))dt\Big)\\
 &=&z^{n+1}{\mathcal L}\Big(\text{p.f.}\int_{0}^{x}\frac{d}{dt}
 \Big(\frac{(x-t)^{n}}{n!}h_{\alpha,k}(t)\Big)dt\\
 && \hspace{3.5cm} + \; \text{p.f.}\int_{0}^{x}
\frac{(x-t)^{n-1}}{(n-1)!}h_{\alpha,k}(t)dt \Big)\\
&=&z^{n+1}{\mathcal L}\Big(\text{p.f.}\int_{0}^{x}\frac{d}{dt}
 \Big(\frac{(x-t)^{n}}{n!}h_{\alpha,k}(t)\Big)dt\Big)+z{\mathcal
 L}(h_{\alpha,k}).
\end{eqnarray*}
The last line results from (4.9). But
$$
\text{p.f.}\int_{0}^{x}\frac{d}{dt}
\Big(\frac{(x-t)^{n}}{n!}h_{\alpha,k}(t)\Big)dt=\sum_{m=0}^{n}\frac{(-1)^{m}x^{n-m}}{(n-m)!m!}
\;
\text{p.f.}\int_{0}^{x}\frac{d}{dt}\Big(t^{m}h_{\alpha,k}(t)\Big)dt.
$$
We get therefore, by (4.7),
\begin{equation*}
\begin{aligned}
\text{p.f.}\int_{0}^{x}\frac{d}{dt}
\Big(\frac{(x-t)^{n}}{n!}h_{\alpha,k}(t)\Big)dt=&\sum_{m=0}^{n}\frac{(-1)^{m}x^{n-m}}{(n-m)!m!}
(h_{\alpha+m,k}^{0}-x^{m}h_{\alpha,k})\\
=& \sum_{m=0}^{n}\frac{(-1)^{m}x^{n-m}}{(n-m)!m!}
(x^{n-m}h_{\alpha+m,k}^{0}-x^{n}h_{\alpha,k})\\ =&
\sum_{m=0}^{n}\frac{(-1)^{m}x^{n-m}}{(n-m)!m!} h_{\alpha+m,k}^{0}.
\end{aligned}
\end{equation*}
The last line results from the fact $\disp\sum_{0\le m\le
n}\frac{(-1)^{m}}{(n-m)!m!}=0$. Hence, by (4.7), we obtain
\begin{equation*}
\begin{aligned}
 \text{p.f.}\int_{0}^{x}\frac{d}{dt}
 \Big(\frac{(x-t)^{n}}{n!}h_{\alpha,k}(t)\Big)dt=&\;
0\;\; \text{if}\;\; \alpha\ne 0,-1,\ldots,-n\;\; \text{or}\;\; k\ne 0\\
=&\;\disp\frac{(-1)^{-\alpha}x^{n+\alpha}}{(n+\alpha)!(-\alpha)!}
\;\;\; \text{otherwise}
\end{aligned}
\end{equation*}
We conclude that for any $\alpha$ and $k$,
\begin{equation}
\begin{aligned}
{\mathcal L}(\frac{d}{dx}h_{\alpha,k})-z{\mathcal
 L}(h_{\alpha,k})=&\;\;
 0\;\;\;\;\;\;\; \text{if}\;\;\;\;\alpha\notin {\bZ}_{<0}\;\; \text{or}\;\; k\ne
0,\\
=&\;\;\disp\frac{(-1)^{-\alpha}z^{-\alpha}}{(-\alpha)!}\;\;\;
\text{otherwise}.
\end{aligned}
\end{equation}
To simplify the notations, we will denote in the sequel $z$ by $x$.
Finally, using the following formula:
$$
\sum_{m=0}^{n}
\frac{(-1)^{m}}{m!(n-m)!}\frac{1}{X+m}=\frac{1}{X(X+1)\ldots(X+n)},
$$
we  deduce from (4.2), (4.3), (4.4) and (4.9) (with $n\ge -\Re
e(\alpha)$) that:\\
 If $\alpha$ is not a negative integer, then
\begin{equation}
\begin{aligned}
{\mathcal L}(h_{\alpha,k})=&\;
\Gamma(\alpha+1)x^{-\alpha-1}\sum_{j=0}^{k}\rho^{(k)}_{\alpha,j}(\ln
x)^j,
                              \;\text{with}\;\;\rho^{(k)}_{\alpha,k}=(-1)^{k}
\;\text{and}\\
&\rho^{(k)}_{\alpha,j}\in
<\Gamma(\alpha),\ldots,\Gamma^{(k)}(\alpha)>_{{\bQ}[\alpha]}\;\;\text{for}\;\;
j=0,\ldots,k-1,
\end{aligned}
\end{equation}
and if $\alpha$ is  a negative integer, then
\begin{equation}
\begin{aligned}
{\mathcal L}(h_{\alpha,k})=&\;
x^{-\alpha-1}\sum_{j=0}^{k+1}\rho^{(k)}_{\alpha,j}(\ln x)^j,
\;\text{with}\;\;\;\rho^{(k)}_{\alpha,k+1}=\frac{(-1)^{\alpha+1}}{(k+1)(-\alpha-1)!},\\
&\;\;\;\rho^{(k)}_{\alpha,k}=\frac{(-1)^{\alpha+1+k}(\alpha+n+2)}{(-\alpha-1)!}+\sum_{m=0\atop
m\ne -\alpha-1}^{n}
\frac{(-1)^{m+k}(\alpha+n+1)}{m!(n-m)!(\alpha+m+1)},\\
&\;\;\;\text{and}\;\;\rho^{(k)}_{\alpha,j}\in
<\Gamma(\alpha),\ldots,\Gamma^{(k)}(\alpha)>_{{\bQ}[\alpha]}\;\;\text{for}\;\;
j=0,\ldots,k-1.
\end{aligned}
\end{equation}
Combining $(4.10)$ and $(4.11)$, we  get  by $x$-adic completion the
following lemma
\begin{lemma} Let $f$ be a finite sum  $\disp\sum_{i}
f_ix^{\alpha_i}(\ln x)^{k_i}x$ where $\alpha_i\in K$, $k_i\in
{\bZ}_{\ge 0}$ and $f_i\in K((x))$. If $f$ is solution of some
operator
 $\phi\in K[x,\disp\frac{d}{dx}]$, then there exists a positive integer $m$
such that $\disp(\frac{d}{dx})^{m}{\mathcal
 F}(\phi)({\mathcal L}(f))=0$. In other word, ${\mathcal L}(f)$ is
a logarithmic solution of $\disp(\frac{d}{dx})^{m}{\mathcal
 F}(\phi)$ at infinity.
\end{lemma}
\subsection{4.2 Arithmetic properties of $\mathcal{L}$}\ \\
In this paragraph, we shall investigate the  relations  between the
radius of convergence of a power series $f\in K((x))$ and those of
the formal factors of $\mathcal{L}(fx^\alpha(\ln x)^k)$ where
$\alpha\in {\bQ}\setminus{\bZ}_{\le 0}$ and $k\in{\bZ}_{\ge 0}$.
\begin{lemma}Let $\alpha\in {\bQ}\setminus{\bZ}_{\le 0}$ and $k\in{\bZ}_{\ge 0}$.
Then, for each $j=0,\ldots,k$, there exist  sequences
$\Big(r_{\alpha+n,j}^{(k,\ell)}\Big)_{n\ge 0}$ of elements of
${\bQ}(\alpha)$, with $\ell=j,\ldots,k$, such that
$$\rho^{(k)}_{\alpha+n,j}=\sum_{\ell=j}^{k}\rho^{(k)}_{\alpha,\ell}\;
r_{\alpha+n,j}^{(k,\ell)}\;\;\;\;
 (n\ge 0),$$
where the $\rho^{(k)}_{\alpha,\ell}$ were defined in (4.12) and
(4.13). Moreover, for any place $v$ of $V_{0}$, these sequences
satisfy
$$\limsup_{n\longrightarrow  \infty}
\Big|r_{\alpha+n,j}^{(k,\ell)}\Big|_{v}^{1/n}\le 1.$$
  \end{lemma}
\begin{proof} We will prove this lemma by downward induction on the
index $j$. In the case $j=k\ge 0$, by (4.12) and (4.13), it suffices
to take $r_{\alpha+n,k}^{(k,k)}=1$ for any $n\in\mathbb{N}$. Suppose
now that the lemma is true for some index $j$ with $1\le j\le k$.
From the formulas (4.12), (4.13) and (4.10), we obtain the
recurrence relation
$$\rho^{(k)}_{\alpha+1,j-1}=-\rho^{(k)}_{\alpha,j-1}+\disp\frac{j}{\alpha+1}\;\rho^{(k)}_{\alpha,j},$$
 and by iteration on $n\ge
1$, we  find
\begin{eqnarray*}
\rho^{(k)}_{\alpha+n,j-1}&=&(-1)^n\rho^{(k)}_{\alpha,j-1}+j\sum_{i=0}^{n-1}\disp\frac{(-1)^{n+i+1}}{\alpha+i+1}\;
\rho^{(k)}_{\alpha+i,j}.\\
&=&(-1)^n\rho^{(k)}_{\alpha,j-1}+j\sum_{i=0}^{n-1}\disp\frac{(-1)^{n+i+1}}{\alpha+i+1}
                              \sum_{\ell=j}^{k}\rho^{(k)}_{\alpha,\ell}\;
 r_{\alpha+i,j}^{(k,\ell)}\\
&=&(-1)^n\rho^{(k)}_{\alpha,j-1}+j\sum_{\ell=j}^{k}\rho^{(k)}_{\alpha,\ell}
\sum_{i=0}^{n-1}\disp\frac{(-1)^{n+i+1}}{\alpha+i+1}\;r_{\alpha+i,j}^{(k,\ell)}.
\end{eqnarray*}
Thus,
\begin{equation}
\begin{aligned}
r_{\alpha+n,j-1}^{(k,\ell)}&=j\sum_{i=0}^{n-1}\disp\frac{(-1)^{n+i+1}}{\alpha+i+1}\;
 r_{\alpha+i,j}^{(k,\ell)}
,\;\;\;\ell=j,\ldots,k\\
r_{\alpha+n,j-1}^{(k,j-1)}&=(-1)^{n},
\end{aligned}
\end{equation}
 we get,
$$\rho^{(k)}_{\alpha+n,j-1}=\sum_{\ell=j-1}^{k}\rho^{(k)}_{\alpha,\ell}\;
 r_{\alpha+n,j-1}^{(k,\ell)}\;\;\;\;\;\;
\mbox{avec}\;\;\;\;\;\; r_{\alpha+n,j-1}^{(k,\ell)}\in
{\bQ}(\alpha).$$ Let $v\in V_{0}$. By induction hypothesis, we have
$\disp\limsup_{n\longrightarrow\infty}
\Big|r_{\alpha+n,j}^{(k,\ell)}\Big|_{v}^{1/n}\le 1$ for
$\ell=j,\ldots,k$. Since $\alpha$ is an element of $K$,  hence
algebraic over ${\bQ}$, it is non-liouville for $p(v)$ and
consequently we have
$\disp\limsup_{n\longrightarrow\infty}\Big|\disp\frac{1}{\alpha+n}\Big|_{v}^{1/n}=1$
  (cf. \cite[VI.1.1]{DGS}). We deduce that
$$\disp\limsup_{n\longrightarrow\infty}\Big(\max_{0\le m\le
n-1}|r_{\alpha+m,j}^{(k,\ell)}\Big|_{v}^{1/n}\Big)\le 1
,\;\;\;\ell=j,\ldots,k$$ and
$$\disp\limsup_{n\longrightarrow\infty} \Big(\max_{0\le m\le
n-1}\Big|\disp\frac{1}{\alpha+m+1}\Big|_{v}^{1/n}\Big)\le 1.$$
Combining these estimations with (4.14) we get for
$\ell=j,\ldots,k$,
$$\disp\limsup_{n\longrightarrow\infty}
\Big|r_{\alpha+n,j-1}^{(k,\ell)}\Big|_{v}^{1/n}\le 1.$$ The case
$\ell=j-1$ is trivial.
\end{proof}
{\bf Notations.}  If $Y\in {\GL}_\mu(K((x)))$, we will denote, for
$s\in {\bZ}$,
\begin{equation*}
\begin{aligned}
 {\mathcal R}_{s}(Y)=&\{y\in K((x))\;|\; r_{v}(y)\ge
R_{v}(Y)\;\pi_{v}^{s},\;\text{
 for almost all}\; v\in V_0\},\\
{\mathcal R}_{s}^\infty(Y)=&\{y(x)\in K((1/x))\;|\; y(1/x)\in
{\mathcal R}_{s}(Y)\}.
\end{aligned}
\end{equation*}
Here, "almost all" means with at most finitely many exceptions. It
is clear that ${\mathcal R}_{s}(Y)$ (resp. ${\mathcal
R}_{s}^\infty(Y)$) is a $K$-subalgebra of $K((x))$ (resp. of
$K((1/x))$). For instance, if $f\in K((x))$, then ${\mathcal
R}_{0}(f)$ denotes the  $K$-algebra  of the power series $y\in
K((x))$ such that $r_{v}(y)\ge r_{v}(f)$ for almost all $v\in V_0$.
\begin{proposition} Let $f\in K[[x]]$ with $f\ne 0$, $
\alpha\in {\bQ}$ and $k\in {\bZ}_{\ge 0}.$ Then there exist power
series $h_{\alpha,k,j}\in {\bC}\otimes_{K}{\mathcal R}_{-1}(f),\;
j=0,\ldots,k$, which satisfy the following conditions
$$
{\mathcal L}\Big(fx^{\alpha}(\ln x)^k\Big)=\left\{
\begin{tabular}{lll}
$x^{-\alpha-1}\Gamma(\alpha)\disp\sum_{j=0}^{k}h_{\alpha,k,j}\Big(\frac{1}{x}\Big)(\ln
x)^j$& if &$\alpha\in {\bQ}\setminus {\bZ}_{<0}$,\\
$\disp\sum_{j=0}^{k+1}h_{\alpha,k,j}\Big(\frac{1}{x}\Big)(\ln x)^j$
& if &$\alpha\in {\bZ}_{<0}$.
\end{tabular}
\right.
$$
with $h_{\alpha,k,k+1}\in K[x]\setminus\{0\}$ and $h_{\alpha,k,k}\in
K[[x]]\setminus\{0\}$ such that
$r_{v}(h_{\alpha,k,k})=r_{v}(f)\pi_{v}^{-1}$ for almost all $v\in
V_0$. In particular, ${\mathcal L}\Big(fx^{\alpha}(\ln x)^k\Big)\ne
0$.
\end{proposition}
\begin{proof} Suppose $\alpha\in {\bQ}\setminus {\bZ}_{<0}$. By
(4.12), we may  write
$${\mathcal
L}\Big(fx^{\alpha}(\ln
x)^k\Big)=x^{-\alpha-1}\Gamma(\alpha)\sum_{j=0}^{k}h_{\alpha,k,j}\Big(\frac{1}{x}\Big)(\ln
x)^j;$$ where
$$h_{\alpha,k,j}=\sum_{n\ge
0}a_{n}\;\frac{\Gamma(\alpha+n+1)}{\Gamma(\alpha)}\;\rho^{(k)}_{\alpha+n,j}x^{n}=
\sum_{n\ge 0}a_{n}(\alpha)_{n+1}\rho^{(k)}_{\alpha+n,j}x^{n}.$$ For
$j=k$, we have $\rho^{(k)}_{\alpha+n,j}=(-1)^{k}$, Thus
$h_{\alpha,k,k}\in K[[x]]$, and since $\alpha\in {\bQ}$, we also
have $\alpha\in{\bZ}_{p(v)}$ for almost all $v\in V_0$. Hence, using
(2.2), we get
$$r_{v}(h_{\alpha,k,k})^{-1}=
\limsup_{n\longrightarrow\infty}|a_{n}(\alpha)_{n+1}|_{v}^{1/n}=
r_{v}(f)^{-1}\pi_{v} \;\;\;\; \text{for almost all}\;\; v\in V_0.$$
For $j=0,\ldots,k-1$, Lemma 4.2 gives
$$
h_{\alpha,k,j}=\sum_{\ell=j}^{k}\rho^{(k)}_{\alpha,\ell}h_{k,j}^{(k,\ell)}
\;\;\;\text{with}\;\;\;h_{\alpha,k,j}^{(k,\ell)}=\sum_{n\ge
0}a_{n}(\alpha)_{n+1}r_{\alpha+n,j}^{(k,\ell)}x^{n}\in K[[x]],
$$ and
$$
r_{v}(h_{\alpha,k,j}^{(k,\ell)})^{-1}\le
\limsup_{n\longrightarrow\infty}|a_{n}(\alpha)_{n+1}|_{v}^{1/n} =
r_{v}(f)^{-1}\pi_{v} \;\;\;\;\text{for almost all}\;\; v\in V_0.
$$
This ends the proof in the case $\alpha\in {\bQ}\setminus
{\bZ}_{<0}$. Now, suppose $\alpha\in {\bZ}_{<0}$, write
$f=\disp\sum_{n\ge 0}a_nx^n$ and put $f_\alpha=\disp\sum_{n\ge
-\alpha}a_nx^n$. Then, by (4.13),
$$\mathcal{L}\Big((f-f_\alpha)x^\alpha (\ln
x)^k\Big)=\sum_{j=0}^{k+1}P_{j}(\frac{1}{x})(\ln x)^j,$$ with
$P_j\in {\bC}[x]$, and $ P_{k+1},P_k \in K[x]\setminus\{0\}.$ On the
other hand, the first assertion, applied to $f_\alpha x^\alpha(\ln
x)^k$, shows that there exist $h_{0,k,j}\in
{\bC}\otimes_{K}{\mathcal R}_{-1}(f),\; j=0,\ldots,k$, such that
$$\mathcal{L}(f_\alpha x^\alpha(\ln
x)^k)=\disp\sum_{j=0}^{k}x^{-1}h_{0,k,j}\Big(\frac{1}{x}\Big)(\ln
x)^j$$ where  $h_{0,k,k}\in K[[x]]$ and
$r_{v}(h_{0,k,k})=r_{v}(f)\pi_{v}$ for almost all $v\in V_0$.
Finally, by linearity of the Laplace transform, it suffices to take
$h_{\alpha,k,j}=xh_{0,k,j}+P_j$ for $j=0,\ldots,k$ and
$h_{\alpha,k,k+1}=P_{k+1}$ to find the second part of the
proposition.
\end{proof}
\section{Solutions of $\psi^*$ and solutions of $\frac{d}{dx}\psi$ }
Let $\psi$  be  a differential operator  of $K[x,d/dx]$ such that
all slopes of $NR(\psi)$ lie in $\{1,0\}$, and $v$ be a fixed
 finite place of $V_0$. Let $K_v$ be the $v$-adic
completion of $K$ and  $\mathrm{\Omega}_{p(v)}$ be a $v$-adic
complete field and algebraically closed containing ${\bC}_{p(v)}$
such that its value group is ${\bR}_{\ge 0}$. We  fix an embedding
$K\hookrightarrow K_v\hookrightarrow {\bC}_{p(v)}\hookrightarrow
\mathrm{\Omega}_{p(v)}$. In this section, we see how we can
determine the nature of solutions of $d/dx.\psi$ at $0$ from those
$\psi^*$ at the same point. For this, we shall begin with the
following key Lemma.
\begin{lemma}
\label{key}Assume $y\in K((x))$, $\alpha\in K\cap {\bZ}_{p(v)}$,
$\delta \in K$ and $k\in {\bZ}_{\ge 0}$. Then the differential
equation $d/dx(z)=yx^\alpha (\ln x)^k\exp(\delta/x)$ has a solution
of the form $\sum_{0\le i\le k+1}y_ix^{\alpha} (\ln
x)^i\exp(\delta/x)$ at $0$, where for $i=0,\ldots,k+1$,  $y_i\in
K((x))$ such that
\begin{equation*}
\begin{aligned}
r_{v}(y_i)\ge& r_v(y)\;\; \text{if} \;\;\delta=0\\
 \ge& \min(|\delta|_v\pi_v^{-1},r_v(y))\;\; otherwise.
\end{aligned}
 \end{equation*}
\end{lemma}
\begin{proof} Let $m=\min(0,\text{ord}_x(y))$ and let write $y=\sum_{n\ge m}a_nx^n\in K((x))$.
Let us consider $z=\disp\sum_{0\le i\le k+1}\sum_{n\ge m}
a_{i,n}x^{n+\alpha} (\ln x)^i\exp(\delta/x)$, where $a_{i,n}\in K$
for all $i=0,\ldots,k+1$ and all $n\ge m$. Then $d/dx(z)=$
$$ \sum_{0\le i\le k}\sum_{n\ge
m+1}\Big((n-1+\alpha)a_{i,n-1}+(i+1)a_{i+1,n-1}-\delta
a_{i,n}\Big)x^{n+\alpha-2}(\ln x)^i\exp(\delta/x)$$
$$+\sum_{n\ge m+1}
\Big((n-1+\alpha)a_{k+1,n-1}-\delta a_{k+1,n}\Big)x^{n+\alpha-2}(\ln
x)^{k+1}\exp(\delta/x)
$$
$$ + \sum_{0\le i\le k+1}-\delta a_{i,m}x^{m+\alpha-2}(\ln
x)^{i}\exp(\delta/x).$$ $z$ is then a solution of the differential
equation $d/dx(z)=yx^\alpha (\ln x)^k\exp(\delta/x)$ if and only if,
the coefficients  $a_{i,n}$ satisfy  the following relations for all
$n\ge m$:
\begin{eqnarray}
&&\delta a_{0,m}=\delta a_{1,m}=\ldots=\delta a_{k+1,m}=0,\\
&&(n+\alpha)a_{k+1,n}-\delta a_{k+1,n+1}=0,\\
&&(n+\alpha+1)a_{k,n+1}+(k+1)a_{k+1,n+1}-\delta a_{k,n+2}=a_n,\\
&&(m+\alpha)a_{k,m}+(k+1)a_{k+1,m}-\delta a_{k,m+1}=0,\\
&&(n+\alpha)a_{i,n}+(i+1)a_{i+1,n}-\delta
a_{i,n+1}=0,\;\;\text{for}\;\; 0\le i<k.
\end{eqnarray}
This means that :\\
\textbf{Case 1.} If $\delta=0$, we have from (5.2) and (5.3),
\begin{equation}
\begin{aligned}
\label{a1} \sum_{n\ge m}a_{k+1,n}x^n=&\; 0\hspace{1.5cm}\text{ if }\;\alpha\;\text{ is  non integer}\;<\;-m,\\
=&\;
a_{k+1,-\alpha}x^{-\alpha}=\disp\frac{a_{-\alpha-1}}{k+1}x^{-\alpha}
\hspace{0.5cm}\text{ otherwise},
\end{aligned}
\end{equation}
and therefore, for all $n\ge m+1$, and all $0\le i \le k$, we get
from (5.3) and (5.5),
\begin{equation}
\begin{aligned}
\label{a2}  a_{1,-\alpha}&=\ldots=a_{k,-\alpha}=0& \;\text{if}\; \alpha\; \text{is an integer}\; \le -m,\\
a_{k,n}&= \frac{a_{n-1}}{n+\alpha}&\;\; \text{for all}\;\; n\ne -\alpha,\\
a_{i,n}&=-\frac{(i+1)a_{i+1,n}}{n+\alpha}=\frac{(-1)^{k-i}k!a_{n-1}}{i!(n+\alpha)^{k-i+1}}&\;\;
\text{for all}\;\; n\ne -\alpha.
\end{aligned}
\end{equation}
Hence, the finite sum $\disp\sum_{0\le i\le k+1}y_ix^{\alpha} (\ln
x)^i$, where the coefficients of the power series $y_i=\sum_{n\ge m}
a_{i,n}x^n$ are defined by \eqref{a1}, \eqref{a2}, where
$a_{0,-\alpha}=0$ if $\alpha$ is an integer $\le -m$, and where
$a_{0,m}=a_{1,m}=\ldots=a_{k+1,m}=0$, is a solution of the equation
$dz/dx=yx^{\alpha} (\ln x)^k$ at $0$. In addition, since $\alpha\in
K\cap {\bZ}_{p(v)}$, it is non-liouville for $p(v)$ and consequently
we have
$\disp\limsup_{n\longrightarrow\infty}\Big|\disp\frac{1}{\alpha+n}\Big|_{v}^{1/n}=1$
  (cf. \cite[VI.1.1]{DGS}). Therefore, by \eqref{a2}, we find for
  $i=0,\ldots,k+1$,
 \begin{eqnarray*}
 \disp\limsup_{n\longrightarrow\infty}|a_{i,n}|_{v}^{1/n}&\le &
    \disp\limsup_{n\longrightarrow\infty}|a_{n-1}|_{v}^{1/n}=
  \limsup_{n\longrightarrow\infty}|a_{n}|_{v}^{1/n}.
\end{eqnarray*}
This implies that $r_v(y_i)\ge r_v(y)$ for $i=0,\ldots,k+1$, and hence the lemma is proved in the case $\delta=0$.\\ \ \\
\textbf{Case 2.} If $\delta\ne 0$, we find, from (5.1) and (5.4),
that $a_{0,m}=a_{1,m}=\ldots=a_{k+1,m}=a_{k,m+1}=0$, and therefore,
by induction on $n\ge m$ and by (5.2), that $\sum_{n\ge m}
a_{k+1,n}x^n=0$. In addition, from (5.3) and (5.5), we get for any
$n \ge m$,
\begin{equation}
\begin{aligned}
\label{b1}
a_{k,n+2}&= \frac{(n+\alpha+1)}{\delta}a_{k,n+1}-\frac{1}{\delta}a_n \\
a_{i,n+1}&=\frac{(n+\alpha)}{\delta}a_{i,n}+\frac{i+1}{\delta}a_{i+1,n}
\;\;\text{for any}\;\; 0\le i<k .
\end{aligned}
\end{equation}
Hence, the finite sum $\disp\sum_{0\le i\le k+1}y_ix^{\alpha} (\ln
x)^i$, where the coefficients of the power series $y_i=\sum_{n\ge m}
a_{i,n}x^n$ are defined recursively  by \eqref{b1} and where
$a_{0,m}=a_{1,m}=\ldots=a_{k,m}=\sum_{n\ge m} a_{k+1,n}x^n=0$, is a
solution of the equation $dz/dx=yx^{\alpha} (\ln x)^k$ at $0$. It
remains to prove that the power series $y_i$ satisfy the condition
of Lemma 5.1.

From \eqref{b1}, we find, for any $n\ge 3$ and any $0\le i< k$, that
\begin{eqnarray*}
a_{k,n+2}&=&
\frac{(n+\alpha+1)(n+\alpha)\ldots(2+\alpha)}{\delta^n}a_{k,2}-\frac{1}{\delta}a_n
-\frac{(n+\alpha+1)}{\delta^2}a_{n-1}\\
&&-\frac{(n+\alpha+1)(n+\alpha)}{\delta^3}a_{n-2}-\ldots-\frac{(n+\alpha+1)(n+\alpha)\ldots(3+\alpha)}{\delta^n}a_{1} \\
\\
a_{i,n+2}&=&\frac{(n+\alpha+1)(n+\alpha)\ldots(2+\alpha)}{\delta^n}a_{i,2}+\frac{i+1}{\delta}a_{i+1,n+1}
\\
&& +\frac{(i+1)(n+\alpha+1)}{\delta^2}a_{i+1,n}+\frac{(i+1)(n+\alpha+1)(n+\alpha)}{\delta^3}a_{i+1,n-1}\\
&&+\ldots+\frac{(i+1)(n+\alpha+1)(n+\alpha)\ldots(3+\alpha)}{\delta^n}a_{i+1,2}.
\end{eqnarray*}
Consequently, if $\alpha$ is non integer $\le -2$, we have for any
$n\ge 1$ and any $0\le i< k$,
\begin{equation}
\begin{aligned}
\label{b} a_{k,n+2}&= \frac{(\alpha+2)_n}{\delta^n}a_{k,2}
-\frac{(\alpha+2)_n}{\delta^{n+1}}\sum_{1\le j\le n}
\frac{\delta^{j}a_j}{(\alpha+2)_j}\\
a_{i,n+2}&=\frac{(\alpha+2)_n}{\delta^n}a_{i,2}
+\frac{(i+1)(\alpha+2)_n}{\delta^{n+2}}\sum_{2\le j\le n+1}
\frac{\delta^{j}a_{i+1,j}}{(\alpha+2)_{j-1}},
\end{aligned}
\end{equation}
and, if $\alpha$ is an integer $\le -2$, we have for $n\ge -\alpha$
and any $0\le i< k$,
\begin{equation}
\begin{aligned}
\label{c} a_{k,n+2}&= \frac{(n+\alpha+1)!}{\delta^n}a_{k,1-\alpha}
-\frac{(\alpha+1+n)!}{\delta^{n+1}}\sum_{-\alpha\le j\le n}
\frac{\delta^{j}a_j}{(\alpha+1+j)!}\\
a_{i,n+2}&=\frac{(n+\alpha+1)!}{\delta^n}a_{i,1-\alpha}
+\frac{(i+1)(\alpha+1+n)!}{\delta^{n+2}}\sum_{-\alpha\le j\le n+1}
\frac{\delta^{j}a_{i+1,j}}{(\alpha+j)!}.
\end{aligned}
\end{equation}
Now,  in the case where $\alpha$ is non integer $\le -2$, we have to
study two cases: \\
\textbf{Case 2.1.a:} If $\alpha$ is non integer $\le -2$, and if
$r_{v}(y)\ge \pi_v^{-1}|\delta|_v$, in other word,
$\disp\limsup_{n\to \infty}|a_n|_v^{1/n}\le \pi_v|\delta|_v^{-1}$,
we have $\disp\limsup_{n\to
\infty}\big(\Big|\frac{\delta^na_n}{(\alpha+2)_n}\Big|_v^{1/n}\Big)\le
1$, since  $\disp\lim_{n\to \infty}|(\alpha+2)_n|_v^{1/n}=\pi_v$.
Then
$$\limsup_{n\to \infty}\Big(\max_{1\le i\le
n}\Big|\frac{\delta^ia_i}{(\alpha+2)_i}\Big|_v^{1/n}\Big)\le 1.$$
This implies that the power series $\disp\sum_{n\ge
2}\Big(\frac{(\alpha+2)_n}{\delta^{n+1}}\sum_{1\le j\le n}
\frac{\delta^{j}a_j}{(\alpha+2)_j}\Big)x^n$ has a radius of
convergence at least $\pi_v^{-1}|\delta|_v$. Thus, by \eqref{b}, we
get $r_{v}(y_k)\ge \pi_v^{-1}|\delta|_v$. Using the same argument,
we prove, by downward induction on the index $i$ and by \eqref{b},
that  $r(y_i)\ge \pi_v^{-1}|\delta|_v$ for
any $0\le i\le k$. This concludes the proof of Lemma 5.1 in the case 2.1.a.\\
\textbf{Case 2.1.b:} If $\alpha$ is non integer $\le -2$, and if
$r_{v}(y)< \pi_v^{-1}|\delta|_v$. We will prove the lemma in this
case by downward induction on the index $i$. First, let $l$ be an
element of $\Omega_{p(v)}$ such that $|l|_v=
\pi_v^{-1}|\delta|_v\disp\limsup_{n\to \infty}|a_n|_v^{1/n}>1$.
Since $\disp\lim_{n\to \infty}|(\alpha+2)_n|_v^{1/n}=\pi_v$, we have
$\disp\limsup_{n\to
\infty}\Big|\frac{\delta^na_n}{l^n(\alpha+2)_n}\Big|_v^{1/n}= 1$,
and hence $$\disp\limsup_{n\to \infty}\Big(\max_{1\le i\le
n}\Big|\frac{\delta^ia_i}{l^i(\alpha+2)_i}\Big|_v^{1/n}\Big)\le 1.$$
Since $|l|_v>1$, we obtain, $\disp\limsup_{n\to
\infty}\Big(\max_{1\le i\le
n}\Big|\frac{\delta^ia_i}{l^n(\alpha+2)_i}\Big|_v^{1/n}\Big)\le 1$ ,
and  $$\disp\limsup_{n\to \infty}\Big(\max_{1\le i\le
n}\Big|\frac{\delta^ia_i}{(\alpha+2)_i}\Big|_v^{1/n}\Big)\le |l|_v=
\pi_v^{-1}|\delta|_v\disp\limsup_{n\to \infty}|a_n|_v^{1/n}.$$ This
shows that the power series $\disp\sum_{n\ge
2}\Big(\frac{(\alpha+2)_n}{\delta^{n+1}}\sum_{1\le j\le n}
\frac{\delta^{j}a_j}{(\alpha+2)_j}\Big)x^n$ has  a radius of
convergence at least $r_{v}(y)$. Thus, by \eqref{b}, we get
$$r_{v}(y_k)\ge \min(\pi_v^{-1}|\delta|_v, r_{v}(y))=r_{v}(y).$$
Suppose now that $r_{v}(y_{i+1})\ge \min(\pi_v^{-1}|\delta|_v,
r_{v}(y))$ for some index $1\le i\le k-1$. \\
If $r_{v}(y_{i+1})< \pi_v^{-1}|\delta|_v$, we find, with the same
argument as above and by \eqref{b}, that  $$r_{v}(y_{i})\ge
\min(\pi_v^{-1}|\delta|_v, r_{v}(y_{i+1}))\ge
\min(\pi_v^{-1}|\delta|_v, r_{v}(y)).$$ If $r_{v}(y_{i+1})\ge
\pi_v^{-1}|\delta|_v$, we get, with the same argument as in case
2.1.a and by \eqref{b},
$$r_{v}(y_{i})\ge \pi_v^{-1}|\delta|_v= \min(\pi_v^{-1}|\delta|_v, r_{v}(y)).$$
This shows that, for all $0\le i\le k$, $r_{v}(y_{i})\ge
\pi_v^{-1}|\delta|_v$. This ends the proof of the lemma in  case 2.1.b .\\
\textbf{Case 2.2:} The case where $\alpha$ is an integer $\le -2$
can be proved with the same arguments employed in cases 2.1.a and
2.1.b, using \eqref{c}, since $\disp\lim_{n\to
\infty}|n!|_v^{1/n}=\pi_v$. This concludes the proof of Lemma 5.1.
\end{proof}
\textbf{Notations:} Let $y_1,\ldots,y_s$ be elements of $K((x))$,
and $\Delta=(\delta_{ij})$ be a $t\times t$ diagonal matrix with
entries in $K$. We denote by $\mathfrak{R}_v(y_1,\ldots,y_s,\Delta)$
the $K$-subalgebra of $K((x))$ consisting of power series $y\in
K((x))$
  satisfying:
\begin{equation*}
\label{ineq}
\begin{aligned}
r_{v}(y)&\ge \min_{1\le h\le s}\{r_v(y_{h})\}\;\; \text{if} \;\;\Delta=0,\\
 &\ge \min(\min_{1\le i\le t}\{|\delta_{ii}|_v\;|\;
  \delta_{ii}\ne 0\}\pi_v^{-1},\min_{1\le h\le s}\{r_v(y_h)\})\;\;
  otherwise.
\end{aligned}
\end{equation*}
Also, we denote by $\mathfrak{R}(y_1,\ldots,y_s,\Delta)$  the
$K$-subalgebra of $K((x))$ consisting of power series $y\in K((x))$
  belonging to $\mathfrak{R}_v(y_1,\ldots,y_s,\Delta)$ for almost all
 $v$ in $V_0$. Again here and in the sequel, "almost all" means with at most finitely many
 exceptions.
\begin{proposition}
\label{sol} Let $\psi\in K[x,d/dx]$ be a differential operator of
rank $\mu$ such that all slopes of $NR(\psi)$ lie in $\{0,1\}$.
Assume that $\psi^*$ has a basis of solutions
 at $0$ with elements in
$$\mathfrak{R}_v(y_1,\ldots,y_s,\Delta)[\ln x,\;
x^{\gamma_{1}},\ldots,x^{\gamma_{s}},
\exp(\delta_{11}/x),\ldots,\exp(\delta_{ss}/x)],$$ where
$y_1,\ldots,y_s$ are elements of $K((x))$, where
$\gamma_1,\ldots,\gamma_s$ are elements of ${\bZ}_{p(v)}\cap K$ and
where $\Delta=(\delta_{ij})$ is a $s\times s$ diagonal matrix with
entries in $K$. Then, $\frac{d}{dx}\psi$ and $(\frac{d}{dx}\psi)^*$
have bases of solutions at $0$ with elements respectively  in $
\mathfrak{R}_v(y_1,\ldots,y_s,\Delta)[\ln x,\;
x^{\pm\gamma_{1}},\ldots,x^{\pm\gamma_{s}},
\exp(\pm\delta_{11}/x),\ldots,\exp(\pm\delta_{ss}/x)]$ and in
$\mathfrak{R}_v(y_1,\ldots,y_s,\Delta)[\ln x,\;
x^{\gamma_{1}},\ldots,x^{\gamma_{s}},
\exp(\delta_{11}/x),\ldots,\exp(\delta_{ss}/x)]. $
\end{proposition}
\begin{proof} Write
$\psi=a_\mu(x)(d/dx)^\mu+a_{\mu-1}(x)(d/dx)^{\mu-1}+\ldots+a_0(x)\in
K[x,d/dx]$. Since all slopes of $NR(\psi)$ lie in $\{0,1\}$, then
$a_\mu$ is a monomial, say $a_\mu=x^\nu$ with $\nu\in {\bZ}_{\ge
0}$, and $\psi$ is regular at infinity. If we denote $\psi_\infty$
the operator obtained from $\psi$ by change of variable $x\to 1/x$,
we find
\begin{eqnarray*}
\psi_\infty&=&
x^{-\nu}(-x^2)^\mu(d/dx)^\mu+\frac{\mu(\mu-1)}{2}(-2x)(-x^2)^{\mu-1}(d/dx)^{\mu-1}\\
&& + (-x^2)^{\mu-1}a_{\mu-1}(1/x)(d/dx)^{\mu-1}+ "\text{terms with
lower degree in}\; d/dx",
\end{eqnarray*}
because for all $a\in K(x)$ and all integer $h\ge 1$, we have
\begin{eqnarray*}
(a.d/dx)^{h}&=&a^{h}(d/dx)^{h}+\frac{h(h-1)}{2}a^{h-1}.d/dx(a).(d/dx)^{h-1}\\
&& +\;\;\text{terms with lower degree in}\; d/dx.
\end{eqnarray*} The
regularity of $\psi_\infty$ at $0$ implies, in particular, that
$\texttt{deg}(a_{\mu-1}(x))\le \nu-1$. In addition, it is easy to
check that
\begin{equation}
\label{lea}
\begin{aligned}
(d/dx.\psi)^*=&\psi^*.d/dx\\
                 =&(-1)^\mu x^\nu(d/dx)^{\mu+1}+((-1)^\mu
                 \nu x^{\nu-1}+(-1)^{\mu-1}a_{\mu-1}(x))(d/dx)^{\mu}\\
                 & +\;"\text{terms with lower
degree in}\; d/dx"\\
                 =&(-1)^\mu x^\nu\Big[(d/dx)^{\mu+1}+(
                 \nu x^{-1}-x^{-\nu}a_{\mu-1}(x))(d/dx)^{\mu}\\
                 & + \;"\text{terms with lower
degree in}\; d/dx"\Big],
\end{aligned}
\end{equation}
because, $(d/dx)^{\mu}.x^\nu=x^\nu(d/dx)^{\mu}+\nu
x^{\nu-1}(d/dx)^{\mu-1}+$ "terms with lower degree in $d/dx$".
 This shows that
\begin{equation}
\label{tr} \texttt{tr}(A_{(d/dx.\psi)^*})=-(\nu
x^{-1}-x^{-\nu}a_{\mu-1}(x))\in \frac{1}{x}K[\frac{1}{x}].
\end{equation}
On the other hand, by hypotheses, $\psi^*$ has a basis of solutions
 $(u_1,\ldots,u_\mu)$ at $0$ such that the $u_i$ are of the form
\begin{equation*}
\begin{aligned}
u_i=&\sum_{\text{finite sum on}\;
j}\widehat{y}_{i_j}x^{\gamma_{i_j}}(\ln
x)^{k_{i_j}}\exp(\delta_{i_j}/x)\\
& \in\; \mathfrak{R}_v(y_1,\ldots,y_s,\Delta)[\ln x,\;
x^{\gamma_{1}},\ldots,x^{\gamma_{s}},
\exp(\delta_{11}/x),\ldots,\exp(\delta_{ss}/x)].
\end{aligned}
\end{equation*}

By Lemma \ref{key}, for $i=1,\ldots,\mu$, the differential equation
$d/dx(z)=u_i$ has a solution of the form
\begin{equation*}
\begin{aligned}
 z_i=&\sum_{\text{finite sum on}\;
j}\Big(\sum_{\text{finite sum on}\;
\ell}\widetilde{y}_{i_{j_\ell}}(\ln
x)^{k_{i_{j_\ell}}}\Big)x^{\gamma_{i_j}}\exp(\delta_{i_j}/x)\\
& \hspace{2cm}\in\; K((x))[\ln x,\;
x^{\gamma_{1}},\ldots,x^{\gamma_{s}},
\exp(\delta_{11}/x),\ldots,\exp(\delta_{ss}/x)]
\end{aligned}
\end{equation*}
such that
\begin{eqnarray*}
 r_{v}(\widetilde{y}_{i_{j_\ell}})&\ge& \{r_v(\widehat{y}_{i_j})\}\;\; \text{if} \;\;\delta_{i_j}=0\\
 &\ge& (|\delta_{i_j}|_v\pi_v^{-1},r_v(\widehat{y}_{i_j}))\;\; otherwise.
 \end{eqnarray*}
Thus, the elements $1,z_1,\ldots,z_\mu$ form a basis of solutions of
$\Big(\frac{d}{dx}\psi\Big)^*=\psi^*\frac{d}{dx}$ at $0$. Moreover,
$1,z_1,\ldots,z_\mu$ lie in
$$\mathfrak{R}_v(y_1,\ldots,y_\mu,\Delta)[\ln
x,\; x^{\gamma_{11}},\ldots,x^{\gamma_{\mu\mu}},
\exp(\delta_{11}/x),\ldots,\exp(\delta_{\mu\mu}/x)].$$

Now, let denote by $W$ the wronskian matrix of $1,z_1,\ldots,z_\mu$.
Thus, the matrices $W$ is solution of $\frac{d}{dx}
X=A_{(\frac{d}{dx}\psi)^*}X$,  and all entries of $W$ lie in
 $$\mathfrak{R}_v(y_1,\ldots,y_\mu,\Delta)[\ln
x,\;x^{\gamma_{1}},\ldots,x^{\gamma_{\mu}},
\exp(\delta_{11}/x),\ldots,\exp(\delta_{\mu\mu}/x)].$$ On the other
hand, We know that $det(W)$ satisfies the differential equation
$d/dx(det(W))=\texttt{tr}(A_{(d/dx.\psi)^*})det(W)$. According to
\eqref{tr}, $det(W)$ is of the form $x^\alpha\exp(P(1/x))$ where
$P\in K[x]$ and where $\alpha\in K$. By definition of $W$, we find
$\alpha \in\; <1,\gamma_{1},\ldots,\gamma_{\mu}>_{\bZ}$ and
$P(x)=\delta x$ for some $\delta \in\;
<\delta_{11},\ldots,\delta_{\mu\mu}>_{\bZ}$. This implies that all
entries of $W^{-1}$ lie in
$$\mathfrak{R}_v(y_1,\ldots,y_\mu,\Delta)[\ln x,\;
x^{\pm\gamma_{1}},\ldots,x^{\pm\gamma_{\mu}},
\exp(\pm\delta_{11}/x),\ldots,\exp(\pm\delta_{\mu\mu}/x)].$$ Hence,
by, Lemma 3.6 and the fact that leading coefficient of
$(d/dx.\psi)^*$ is monomial (see \eqref{lea}), the differential
operator $\frac{d}{dx}\psi$ has a basis of solutions at $0$, with
element in
$$\mathfrak{R}_v(y_1,\ldots,y_\mu,\Delta)[\ln
x,\; x^{\pm\gamma_{1}},\ldots,x^{\pm\gamma_{\mu}},
\exp(\pm\delta_{11}/x),\ldots,\exp(\pm\delta_{\mu\mu}/x)].$$ This
concludes the proof of Proposition 5.2.
\end{proof}
\begin{corollary}
Under hypotheses of Proposition 5.2, for all positive integer $m\ge
1$, the differential operators $(\frac{d}{dx})^m\psi$ and
$((\frac{d}{dx})^m\psi)^*$ have bases of solutions at $0$ with
elements respectively  in
$$
\mathfrak{R}_v(y_1,\ldots,y_\mu,\Delta)[\ln x,\;
x^{\pm\gamma_{1}},\ldots,x^{\pm\gamma_{\mu}},
\exp(\pm\delta_{11}/x),\ldots,\exp(\pm\delta_{\mu\mu}/x)]$$ and
$\;\;\mathfrak{R}_v(y_1,\ldots,y_\mu,\Delta)[\ln x,\;
x^{\gamma_{1}},\ldots,x^{\gamma_{\mu}},
\exp(\delta_{11}/x),\ldots,\exp(\delta_{\mu\mu}/x)].$
\end{corollary}
\begin{proof}First, the differential operator $(\frac{d}{dx})^m\psi$
has the same leading coefficient as $\psi$ which is a monomial. In
addition, the properties of Newton polygon (\cite[III.1]{Ma}) lead
to
$$\{\text{slopes of }\; N\Big((\frac{d}{dx})^m\psi\Big)\}=\{\text{slopes of
}\; N(\frac{d}{dx})\}\cup\{\text{slopes of }\;
N(\psi)\}\in\{0,1\},$$ and
\begin{equation*}
\begin{aligned}
\{\text{slopes of }\;
N\Big(((\frac{d}{dx})^m\psi)_\infty\Big)\}=&\{\text{slopes of }\;
N((\frac{d}{dx})_\infty)\}\cup\{\text{slopes of }\;
N(\psi_\infty)\}\\
=&\{0\}.
\end{aligned}
\end{equation*}
Thus, the slopes of $NR((\frac{d}{dx})^m\psi)$ lie in $\{0,1\}$ for
all integer $m\ge 1$. Hence, the corollary can be proved by
induction on $m$, using Proposition 5.2.
\end{proof}
Let $f_E$ denote the Euler series : $\sum_{n\ge 0}(-1)^{n}n!x^n$.
With the notations of \S4.2, we obtain:
\begin{corollary}
Let $\psi\in K[x,d/dx]$ be a differential operator of rank $\mu$
such that all slopes of $NR(\psi)$ lie in $\{0,1\}$. Assume that the
differential system $dX/dx=A_{\psi}X$ has solution at $0$ of the
form $Y(x)x^\Gamma\exp(\Delta/x)$,
 where the $Y(x)$ is a $\mu\times\mu$ invertible matrix with entries in $K((x))$, where
 $\Gamma$ is a $\mu\times\mu$ matrix with entries in $K$ and
 eigenvalues $\gamma_1,\ldots,\gamma_\mu$ in ${\bQ}$, and where
 $\Delta=(\delta_{ij})$ is a $\mu\times\mu$ diagonal
 matrix with entries in $K$. Then, for all positive integer $m\ge
1$, the differential operators $(\frac{d}{dx})^m\psi$ and
$((\frac{d}{dx})^m\psi)^*$ have bases of solutions at $0$ with
elements respectively  in
$$
\Big(\mathcal{R}_0(Y)\cap\mathcal{R}_0(f_E)\Big)[\ln x,\;
x^{\pm\gamma_{1}},\ldots,x^{\pm\gamma_{\mu}},
\exp(\pm\delta_{11}/x),\ldots,\exp(\pm\delta_{\mu\mu}/x)]$$ and
$\;\;\Big(\mathcal{R}_0(Y)\cap\mathcal{R}_0(f_E)\Big)[\ln x,\;
x^{-\gamma_{1}},\ldots,x^{-\gamma_{\mu}},
\exp(-\delta_{11}/x),\ldots,\exp(-\delta_{\mu\mu}/x)].$
\end{corollary}
\begin{proof}Since all slopes of $NR(\psi)$ lie in $\{0,1\}$, the leading
coefficient $a_\mu$ of $\psi$ is a monomial. Let
$\widetilde{Y}_{\mu}(x)=(\widetilde{y}_1(x),\ldots,\widetilde{y}_\mu(x))\in
{\M}_{\mu\times 1}(K((x)))$ denotes the $\mu$-th row of the matrix
$a_{\mu}^{-1}(^TY(x)^{-1})$. By Lemma 3.7, the elements of
$\widetilde{Y}_{\mu}(x)(x)^{(-^T\Gamma)}\exp(-^T\Delta/x)$ form a
basis of solutions of $\psi^*$ at $0$. According to Corollary 3.3,
the elements of this basis lie, for all $v\in V_0$, in
$$\mathfrak{R}_v(\widetilde{y}_1,\ldots,\widetilde{y}_\mu,\Delta)[\ln
x,\; x^{-\gamma_{1}},\ldots,x^{-\gamma_{\mu}},
\exp(-\delta_{11}/x),\ldots,\exp(-\delta_{\mu\mu}/x)].$$ In
addition, for almost all $v\in V_0$, the eigenvalues of $\Gamma$ lie
in ${\bZ}_{p(v)}$. Hence, by Corollary 5.3,  the differential
operators $(\frac{d}{dx})^m\psi$ and $((\frac{d}{dx})^m\psi)^*$ have
bases of solutions at $0$ with elements respectively  in
$$
\mathfrak{R}(\widetilde{y}_1,\ldots,\widetilde{y}_\mu,\Delta)[\ln
x,\; x^{\pm\gamma_{1}},\ldots,x^{\pm\gamma_{\mu}},
\exp(\pm\delta_{11}/x),\ldots,\exp(\pm\delta_{\mu\mu}/x)]
$$
and
$\mathfrak{R}(\widetilde{y}_1,\ldots,\widetilde{y}_\mu,\Delta)[\ln
x,\; x^{-\gamma_{1}},\ldots,x^{-\gamma_{\mu}},
\exp(-\delta_{11}/x),\ldots,\exp(-\delta_{\mu\mu}/x)].$\\ The
corollary results therefore from the following observation :
$$\mathfrak{R}(\widetilde{y}_1,\ldots,\widetilde{y}_\mu,\Delta)\subseteq
\mathcal{R}_0(\widetilde{y}_1)\cap\ldots\mathcal{R}_0(\widetilde{y}_\mu)\cap\mathcal{R}_0(f_E)
\subseteq\mathcal{R}_0(Y)\cap\mathcal{R}_0(f_E).$$
\end{proof}
\section{Sufficient conditions}
Let  $\overline{\mathcal{F}}$ denote the inverse of $\mathcal{F}$,
that is the $K$-automorphism of $K[x,d/dx]$ satisfying
$\overline{\mathcal{F}}(x)=-d/dx$ and
$\overline{\mathcal{F}}(d/dx)=x$. In this section, we will prove
that the conditions in Theorem 3.1 are sufficient:
\begin{theorem} Let $\psi\in K[x,d/dx]$ be an operator  of rank $\mu$
satisfying the following conditions:
\begin{itemize}
\item[(1)] the coefficients of $\psi$ are not all in K;
\item[(2)] the slopes of $NR(\psi)$ lie  in $\{-1,0\}$, \item[(3)]
the differential system $d/dxZ=A_{\psi}Z$  has a solution of the
from
$$Y_\psi(\frac{1}{x})(\frac{1}{x})^\Gamma\exp(-\Delta x),$$
where $Y_\psi(x)$ is a $\mu\times\mu$ invertible matrix  with
entries in $K((x))$ such that $\prod_{v\in V_{0}} \min(
R_{v}(Y_\psi)\pi_{v},1)\ne 0,$ where $\Gamma$ is  a $\mu\times\mu$
matrix with entries in $K$ and with eigenvalues
$\gamma_1,\ldots,\gamma_\mu$ in $\mathbb{Q}$, and where
$\Delta=(\delta_{ij})$ is a
 diagonal $\mu\times\mu$ matrix with entries in $K$ which commutes
with $\Gamma$.
\end{itemize}
Then $\psi$ is an $E$-operator.
\end{theorem}
Note that the condition (1) means that the differential operator
$\phi:=\overline{\mathcal{F}}(\psi)$ is not a polynomial.
\begin{lemma}Under the hypotheses of Theorem 5.2, the differential
operator $\phi:=\overline{\mathcal{F}}(\psi)$ has a basis of
solutions at $0$ of the form $(f_1,\dots,f_\nu)x^{C}$ where
$f_1,\dots,f_\nu$ are power series of $K[[x]]$ such that
$\prod_{v\in V_{0}}\min(r_{v}(f_i),1)\ne 0$ for $i=1,\ldots,\nu$ and
where $C$ is a $\nu\times\nu$ upper triangular matrix with entries
in ${\bQ}$.
\end{lemma}
\begin{proof} By \S2.2, $\phi$ is regular at $0$ and admits a  basis of
solutions at $0$ of the form
\begin{eqnarray*}
(\zeta_{1},\zeta_{2},\ldots,\zeta_{\nu})&:=&
(f_{1},f_{2},\ldots,f_{\nu})x^{C}
\end{eqnarray*}
such that \\
1)$f_{1},\ldots,f_{\nu}\in K[[x]]$,\\
 2) $C=D+N$  is an $\nu\times\nu$ matrix where $D$ is a diagonal matrix
 whose diagonal entries
  $D_{ii}:=\alpha_{i}\in K$ are the exponents of $\phi$ at
$0$, and
  $N=(N_{ij})$ is an upper triangular nilpotent  matrix with
entries in  ${\bQ}$ such that $DN=ND$.

Since $ x^{C}= x^{D+N}=x^{D}\sum_{k\ge 0}\disp\frac{N^{k}}{k!}(\ln
x)^k =x^{D}+x^{D}\sum_{k=1}^{\nu}\disp\frac{N^{k}}{k!}(\ln x)^k, $
we obtain
\begin{eqnarray} \zeta_{1}&=&f_{1}x^{\alpha_{1}},
\end{eqnarray} and for $1<i\le \nu$
\begin{eqnarray} \zeta_{i}&=&f_{i}x^{\alpha_{i}}+\sum_{j=1}^{
i-1}f_{j}x^{\alpha_{j}}\sum_{k=1}^{\nu}\disp\frac{(N^{k})_{ji}}{k!}(\ln
x)^k,
\end{eqnarray}
since $(N^{k})_{ji}=0$  for  $j\ge i$.

In addition, by Lemma 4.1, there exists a positive integer $m$ such
that $\disp\Big(\frac{d}{dx}\Big)^{m}\psi$ annihilates  ${\mathcal
L}(\zeta_{i})$ for $i=1,\ldots,\nu$. We then define,
$\Psi=\disp\Big(\frac{d}{dx}\Big)^{m}\psi$. Applying Corollary 5.4
to $\psi$ at infinity, we find that, $\Psi$ has a basis  of
solutions $ \xi_1,\ldots,\xi_{\mu+m}$ at infinity with elements in
$$
\Big(\mathcal{R}_0^\infty(Y_\psi)\cap\mathcal{R}_0^\infty(f_E)\Big)[\ln
x,\; x^{\pm\gamma_{1}},\ldots,x^{\pm\gamma_{\mu}},
\exp(\pm\delta_{11}x),\ldots,\exp(\pm\delta_{\mu\mu}x)].$$

Now, let  ${\mathcal A}_{0}$ denote the set
$${\bC}\otimes_{K}\Big(\mathcal{R}_0^\infty(Y_\psi)\cap\mathcal{R}_0^\infty(f_E)\Big)[\ln x,x^{\pm\gamma_{1}},\ldots,x^{\pm\gamma_{\mu}},
\exp(\pm\delta_{11}x),\ldots,\exp(\pm\delta_{\mu\mu}x)].
$$
Therefore, we have for all $1\le i\le \nu$, ${\mathcal
L}(\zeta_{i})\in\; {\mathcal A}_{0}.$ By induction on $i$, we deduce
from (6.1) and (6.2) that
\begin{eqnarray}\mathcal{L}\Big(f_{i}x^{\alpha_{i}}\Big)\in
{\mathcal A} _{i-1}\;\;\;\;\;(i=1,\ldots,\nu)
\end{eqnarray}
where ${\mathcal A}_{1},\ldots,{\mathcal A}_{\nu-1}$ are the $
{\bC}[\ln x]$-modules of finite type defined recursively by
\begin{eqnarray*}
{\mathcal A}_{i}&=&{\mathcal A}_{i-1} +
\Big<\mathcal{L}(f_{i}x^{\alpha_{i}}(\ln x)^j);\; 0\le j\le
\nu\Big>_{\mathbb{C}[\ln x]}.
\end{eqnarray*} This shows, by iteration on $i$ and by (4.12) and (4.13), that the exponents
 $\alpha_{i}$ are rational numbers. Thus, by
Proposition 4.3, the Laplace
 transform of $f_{i}x^{\alpha_{i}}(\ln x)^k$ (for $i=1,\ldots,\nu$
and $k\in{\bZ}_{\ge 0}$) can be written as
\begin{equation}
\begin{aligned}
{\mathcal L}\Big(f_ix^{\alpha_i}(\ln x)^k\Big)&=
x^{-\alpha_i-1}\Gamma(\alpha)\disp\sum_{j=0}^{k}h_{i,k,j}\Big(\frac{1}{x}\Big)(\ln
x)^j\;\; \text{if}\;\;\alpha_i\in {\bQ}\setminus {\bZ}_{<0},\\
&=\disp\sum_{j=0}^{k+1}h_{i,k,j}\Big(\frac{1}{x}\Big)(\ln
x)^j\;\;\text{if}\;\;\alpha_i\in {\bZ}_{<0}.
\end{aligned}
\end{equation}
where $h_{i,k,j}\in {\bC}\otimes_{K}{\mathcal R}_{-1}(f_i),\;
j=0,\ldots,k$, $h_{i,k,k+1}\in K[x]\setminus\{0\}$ and $h_{i,k,k}\in
K[[x]]\setminus\{0\}$ such that
$r_{v}(h_{i,k,k})=r_{v}(f_i)\pi_{v}^{-1}$ for almost all
$v\in V_0$.\\
To conclude, it  suffices to prove, by induction on $i$, that
\begin{eqnarray} f_{i}\in {\mathcal
R}_{1}(Y_{\psi})\cap {\mathcal R}_{1}(f_E),\;\;\;\;
(i=1,\ldots,\nu).
\end{eqnarray}

Combining  (6.3) with (6.4) for  $i=1$ and $k=0$, we find that
$\alpha_{1}\in\{\pm\gamma_{j}+m\;|\; m\in{\bZ},\; j=1,\ldots,\mu\}$
and that $h_{1,0,0} \in {\mathcal R}_{0}(Y_{\psi})\cap {\mathcal
R}_{0}(f_E).$ Thus $ f_{1}\in {\mathcal R}_{1}(Y_{\psi})\cap
{\mathcal R}_{1}(f_E)$ and for any $0\le j\le k$, we have
\begin{eqnarray*} h_{1,k,j}\in {\bC}\otimes_{K}({\mathcal
R}_{0}(Y_{\psi})\cap {\mathcal R}_{0}(f_E)),
\end{eqnarray*}
and hence, for any $k\ge 0$,
\begin{eqnarray*}
\mathcal{L}\Big(f_{1}x^{\alpha_{1}}(\ln x)^k\Big)\in
x^{-\alpha_{1}}{\bC}\otimes_{K}({\mathcal R}_{0}(Y_{\psi})\cap
{\mathcal R}_{0}(f_E))[\ln x].
\end{eqnarray*}
This implies  $${\mathcal A}_{1}\;\subseteq\;{\mathcal A}_{0}.
$$ Suppose now that, for some
integer $\tau$ with $1\le\tau-1<\nu$, we have $f_{i}\in {\mathcal
R}_{1}(Y_{\psi})\cap {\mathcal R}_{1}(f_E),$ and
$\alpha_{i}\in\{\pm\gamma_{j}+m\;|\; m\in{\bZ},\; j=1,\ldots,\mu\},$
for $ i=1,\ldots,\tau-1.$
 Then, by (6.4),
$$h_{i,k,j}\in {\bC}\otimes_{K}{\mathcal
R}_{0}(Y_{\psi})\cap {\mathcal
R}_{0}(f_E))\;\;\;\;\text{for}\;\;\;1\le
 i\le\tau-1,\;\text{and }\; 0\le j\le k.$$
This implies $
 {\mathcal A}_{\tau-1}\;\subseteq\;{\mathcal A}_{0}.
$ In particular, by (6.3), we get,
$\mathcal{L}\Big(f_{\tau}x^{\alpha_{\tau}}\Big)\;\in {\mathcal
A}_{0}.$ Therefore, by (6.4), we find
$\alpha_{\tau}\in\{\pm\gamma_{j}+m\;|\; m\in{\bZ},\;
j=1,\ldots,\mu\}$ and $h_{\tau,0,0} \in {\mathcal
R}_{0}(Y_{\psi})\cap {\mathcal R}_{0}(f_E),$ and consequently $
f_{\tau}\in {\mathcal R}_{1}(Y_{\psi})\cap {\mathcal R}_{1}(f_E).$
This prove that $f_{i} \in {\mathcal R}_{1}(Y_{\psi})\cap {\mathcal
R}_{1}(f_E)$ and $\alpha_{i}\in\{\pm\gamma_j+m\;|\; m\in{\bZ},\;
j=1,\ldots,\mu\}$ for $i=1,\ldots,\nu$.  On the other hand, by
Corollary 3.4, the power series $f_{i}$ are entries  of the inverse
of a reduction matrix of $A_\phi$. Then, by Proposition 2.1, they
satisfy $r_v(f_i)\ne 0$ for any $v\in V_0$. Combining this with the
fact that $f_{i} \in {\mathcal R}_{1}(Y_{\psi})\cap {\mathcal
R}_{1}(f_E)$ for $i=1,\ldots,\nu$, we get $$\disp\prod_{v\in
V_{0}}\min(r_{v}(f_i),1)\ne 0\;\;\;\text{ for}\;\;i=1,\ldots,n.$$
The lemma follows therefore since $\alpha_1,\ldots,\alpha_\nu\in
{\bQ}$.
\end{proof}
\begin{proof}[Proof of Theorem 5.2] First, by \S2.3, the differential operator $\phi^{*}$ is regular at $0$.
 In addition, by Corollary 3.7, the differential system
$dX/dx=A_{\overline{\psi^*}}X$ has a solution of the form
$\widetilde{Y}(\frac{1}{x})(\frac{1}{x})^{-^T\Gamma}\exp(\Delta x)$,
where $\widetilde{Y}(x)\in {\GL}_{\mu}(K((x)))$ such that
$\disp\prod_{v\in V_{0}}\min(r_{v}(\widetilde{Y}),1)\ne 0$.
Moreover, we have ${\mathcal F}(\phi^{*})=({\mathcal
F}\overline{\phi})^{*}=\overline{\psi^{*}}$ (cf. \cite[V.3.6]{Ma})).
Then, by the same proof as in Lemma 6.2, we find that $\phi^{*}$ has
also a basis of solutions at infinity of the form
$(z_1,\dots,z_\nu)x^{\Lambda}$ where $z_1,\dots,z_\nu$ are power
series of $K[[x]]$ such that $\disp\prod_{v\in
V_{0}}\min(r_{v}(z_i),1)\ne 0$  for $i=1,\ldots,\nu$ and where
$\Lambda$ is a $\nu\times\nu$ upper triangular matrix with entries
in  ${\bQ}$. Combining this with Lemma 6.2 and Lemma 3.4, we find
that the differential system $dX/dx=A_\phi X$ has a solution at $0$
of the form  $Y_\phi(x)x^C$ where $Y(x)\in {GL}_{\nu}(K((x)))$  such
that $\disp\prod_{v\in V_{0}}\min(r_{v}(Y_\phi),1)\ne 0$ and where
$C\in {M}_\nu({\bQ})$ is an upper triangular matrix (see proof of
Theorem 3.8). Hence, $\phi$ is a $G$-operator and consequently
$\psi$ is an $E$-operator.
\end{proof}

\end{document}